\documentstyle[a4,12pt, eucal, amstex, amssymb, latexsym]{amsart}
\newtheorem{theo}{Theorem}[section]
\newtheorem{lemm}[theo]{Lemma}
\newtheorem{cor}[theo]{Corollary}
\newtheorem{prop}[theo]{Proposition}

\newtheorem{rem}{Remark}[section]
\newcommand{\Der}{{\rm Der\:}}
\newcommand{\ad}{{\rm ad\:}}

\newcommand{\Z}{{\Bbb Z}}
\newcommand{\F}{{\Bbb F}}
\newcommand{\OO}{{\mathcal O}}
\newcommand{\rad}{{\rm rad}}
\newcommand{\un}{\underline}
\newcommand{\ot}{\otimes}
\begin{document}

\title{Simple Lie algebras of small characteristic V. The non-Melikian case}
\author{\sc Alexander Premet and Helmut Strade}
\address
{School of Mathematics, The University of Manchester, Oxford Road,
M13 9PL, United Kingdom} \email{sashap@@maths.man.ac.uk}
\address{
Fachbereich Mathematik, Universit{\"a}t Hamburg, Bundesstrasse 55,
20146 Hamburg, Germany}
\email{strade@@math.uni-hamburg.de}
\thanks{
Primary 17B20, 17B50}
\begin{abstract}
\noindent Let $L$ be a finite-dimensional simple Lie algebra over
an algebraically closed field $F$ of characteristic $p>3$. We
prove in this paper that if for every torus $T$ of maximal
dimension in the $p$-envelope of $\ad L$ in $\Der L$ the
centralizer of $T$ in $\ad L$ acts triangulably on $L$, then $L$
is either classical or of Cartan type. As a consequence we obtain
that {\it any} finite-dimensional simple Lie algebra over an
algebraically closed field of characteristic $p>5$ is either
classical or of Cartan type. This settles the last remaining case
of the generalized Kostrikin-Shafarevich conjecture (the case
where $p=7$).
\end{abstract}
\maketitle

\section{Introduction and preliminaries}
Let $F$ be an algebraically closed field of characteristic $p>3$.
In this note we go through the relevant parts of the second
author's classification of finite-dimensional simple Lie algebras
of characteristic $p>7$ to check whether the results there need
additions, modifications or supplementary proofs in order to apply
in the present case where $p>3$. It turned out that only few
changes are necessary.

In what follows $L$ will always denote a finite-dimensional simple
Lie algebra over $F$. We identify $L$ with the subalgebra $\ad L$
of the derivation algebra $\Der L$ and let $L_p$ be the
$p$-envelope of $L$ in the restricted Lie algebra $\Der L$. We
denote by $T$ a torus of maximal dimension in $L_p$ and set
$$\widetilde{H}:=\,{\mathfrak c}_{L_p}(T)=\{x\in L_p\,|\,\,[t,x]=0\
\,\mbox{for all } \,t\in T\},\quad H :=\,{\mathfrak
c}_L(T)\,=\,\widetilde{H}\cap L.$$ A torus $T$ is called {\it
standard} if $H^{(1)}$ consists of nilpotent derivations of $L$.
We denote by $\Gamma(L,T)$ the set of roots of $L$ relative to $T$
(roots are {\it nonzero} linear functions $\gamma\in T^*$ such
that $L_\gamma:=\{x\in L\,|\,\,[t,x]=\gamma(t)x\ \,\mbox{for all }
t\in T\}$ is nonzero). We have root space decompositions
$$L\,=\,H\oplus\bigoplus_{\gamma\in\Gamma(L,T)}L_\gamma,\quad\
\,\,\,
L_p\,=\,\widetilde{H}\oplus\bigoplus_{\gamma\in\Gamma(L,T)}L_\gamma.$$

By \cite[Corollary 3.7]{PS4}, only four types of roots can occur
in simple Lie algebras of characteristic $p>3$: solvable,
classical, Witt, and Hamiltonian roots. In other words, for any
$\gamma\in\Gamma(L,T)$ the semisimple quotient
$L[\gamma]=L(\gamma)/\rad\, L(\gamma)$ of the $1$-section
$L(\gamma):=H\oplus\bigoplus_{i\in\F_p}L_{i\gamma}$ is either
$(0)$ or ${\mathfrak sl}(2)$ or the Witt algebra or contains an
isomorphic copy of the Hamiltonian algebra $H(2;\un{1})^{(2)}$ as
an ideal of codimension $\le 1$. The main result of \cite{PS4}
states that if $L_p$ contains a torus $T'$ of maximal dimension
such that all roots in $\Gamma(L,T')$ are solvable or classical,
then $L$ is either a classical Lie algebra or a filtered Lie
algebra of type $S$ or $H$; see \cite[Theorems C and D]{PS4}.

\smallskip

In this note we impose the following two assumptions on $L$:

\smallskip

\begin{itemize}
\item [$\bullet$] all tori of maximal dimension in $L_p$ are
standard;

\smallskip

\item [$\bullet$] the set of roots of any torus of maximal
dimension in $L_p$  contains either a Witt root ar a Hamiltonian
root.
\end{itemize}

\begin{theo}
If a finite dimensional simple Lie algebra $L$ over $F$ satisfies
the above assumptions, then $L$ is isomorphic to a Lie algebra of
Cartan type.
\end{theo}
Combining Theorem~1.1 with \cite[Theorems C and D]{PS4} we derive:
\begin{theo}
Let $L$ be a finite dimensional simple Lie algebra over an
algebraically closed field of characteristic $p>3$ such that all
tori of maximal dimension in  $L_p$ are standard. Then $L$ is either
classical or of Cartan type.
\end{theo}
Due to \cite{Wil77, St89a, P94} the assumption on tori in
Theorem~1.2 is fulfilled automatically when $p>5$. In this case
Theorem~1.2 can be rephrased as follows:
\begin{theo}
Any finite dimensional simple Lie algebra over an algebraically
closed field of characteristic $p>5$ is either classical or of
Cartan type.
\end{theo}
Theorem~1.3 settles the last remaining case $p=7$ of the
Kostrikin--Shafarevich conjecture on the structure of
finite-dimensional restricted simple Lie algebras over
algebraically closed fields of characteristic $p>5$; see
\cite{KS66}. In the early 80s, G.M.~Melikian discovered a
restricted simple Lie algebra of characteristic $5$ which was
neither classical nor of Cartan  type, thereby showing that the
restriction on $p$ in the Kostrikin--Shafarevich conjecture was
necessary. In 1984, R.E.~Block and R.L.~Wilson succeeded to prove
the Kostrikin--Shafarevich conjecture for algebraically closed
fields of characteristic $p>7$; see \cite{BW}.

As far as the general classification problem for $p>3$ is concerned,
Theorem~1.2 leaves open the case where $p=5$ and $L_p$ contains
nonstandard tori of maximal dimension. This isolated case will be
treated in \cite{PS6}, the last paper of the series.

\section{Two-sections of $L$}
In the next two sections our standing hypothesis is that $L$ is a
finite dimensional simple Lie algebra such that all tori of maximal
dimension in $L_p$ are standard. The second assumption of Sect.~1
will come into force in Sect.~4. We retain the notation introduced
in \cite{PS1}, \cite{PS2}, \cite{PS3}, \cite{PS4} with the following
two exceptions: to match the notation of \cite{St04} we will denote
the divided power algebra $A(m;\un{n})$ by ${\mathcal O}(m;\un{n})$
and the Melikian algebra ${\frak g}(m,n)$ by ${\cal M}(m,n)$.

Our first result extends \cite[Theorems 3.1, 1.7, 1.8]{St89b} and
\cite[Corollary 1.9]{BOSt} which hold for $p>7$.
\begin{theo}
Let $T$ be any torus of maximal dimension in $L_p$.
\begin{itemize}
\item[(i)] The subalgebra $\widetilde{H}={\mathfrak c}_{L_p}(T)$
acts triangulably on $L$.

\smallskip

\item[(ii)] For every $\gamma\in\Gamma(L,T)$ the radical $\rad\,
L(\gamma)$ is $T$-invariant and the factor algebra
$L[\gamma]=L(\gamma)/\rad\, L(\gamma)$ is either zero or
isomorphic to one of ${\mathfrak sl}(2)$, $W(1;\underline{1})$,
$H(2;\underline{1})^{(2)},$ $H(2;\underline{1})^{(1)}$.
\end{itemize}
\end{theo}
\pf Since all tori of maximal dimension in $L_p$ are assumed to be
standard, the first statement is nothing but \cite[Theorem
3.12]{PS4}, while the second statement is immediate from
\cite[Corollary 3.7]{PS4}. \qed

Given a filtered Cartan type Lie algebra $\mathfrak g$ (not
necessarily simple) we denote by ${\mathfrak g}_{(0)}$ the standard
maximal subalgebra of $\mathfrak g$. When ${\mathfrak g}=
W(1;\un{1})$, we have $\dim ({\mathfrak g}/{\mathfrak g}_{(0)})=1,$
while when ${\mathfrak g}\cong H(2;\un{1})^{(\epsilon)}$ with
$\epsilon=1,2,$ we have $\dim ({\mathfrak g}/{\mathfrak
g}_{(0)})=2$. Theorem 2.1 shows that every $1$-section $L(\gamma)$
with $\gamma\in\Gamma(L,T)$ contains a distinguished subalgebra
$Q(\gamma)$ such that $Q(\gamma)/\rad\, Q(\gamma)$ is either zero or
isomorphic to ${\mathfrak sl}(2)$. More precisely, the following
holds:
\begin{enumerate}
\item[(a)] $Q(\gamma)\,=\,L(\gamma)$ if $L(\gamma)$ is solvable or
$L[\gamma]\cong {\mathfrak sl}(2)$; \item[(b)]
$(Q(\gamma)+\rad\,L(\gamma))/\rad\,L(\gamma)\,=\,L[\gamma]_{(0)}$ if
$L[\gamma]$ is of Cartan type.
\end{enumerate}
This is analogous to \cite[Proposition 1.9]{St89b} and
\cite[Proposition 1.11]{BOSt}.

Recall that a root $\gamma\in\Gamma(L,T)$ is called {\it proper}, if
$Q(\gamma)$ is $T$-invariant, and {\it improper} otherwise; see
\cite{PS4}. This definition differs from that introduced by
Block--Wilson. However, it agrees with the Block--Wilson definition
when $p>7$ and reflects better the desired properties of $\gamma$
when $p\in\{5,7\}$ (the formal extension of the Block--Wilson
definition to the case where $p=5$ would imply that all Hamiltonian
roots are improper, which is not what we want).

\medskip

The following result is very important for the classification:
\begin{theo}
Let $L(\alpha,\beta)$ be any $2$-section of $L$ relative to a torus
$T$ of maximal dimension in $L_p$, and $I(\alpha,\beta)$ the maximal
solvable ideal of $T+L(\alpha,\beta)$. Let $\psi$ denote the
canonical homomorphism $T+L(\alpha,\beta)\rightarrow
\big(T+L(\alpha,\beta)/I(\alpha\,\beta)\big)$, and put
$K:=\psi(L(\alpha,\beta))$ and $\overline{T}:=\psi(T)$. Then one of
the following holds:
\begin{itemize}
\item[(1)] $L(\alpha,\beta)$ is solvable and $K=(0)$;

\smallskip

 \item[(2)] $S=\overline{T}+K$ where $S$ is one of ${\mathfrak
sl}(2)$, $W(1;\un{1})$, $H(2;\un{1})^{(2)}$ or else
$S=H(2;\un{1})^{(2)}$ and $S\subset \overline{T}+K\subset
H(2;\un{1})^{(1)}$. Moreover, there exists a root
$\mu\in\F_p\alpha+\F_p\beta$ such that $K=\psi(L(\mu))$;

\smallskip

\item[(3)] $S_1\oplus S_2\subset \overline{T}+K\subset(\Der
S_1)^{(1)}\oplus(\Der S_2)^{(1)}$ where $S_i$ is one of
${\mathfrak sl}(2)$, $W(1;\un{1})$, $H(2;\un{1})^{(2)}$ for
$i=1,2$;

\smallskip

\item[(4)] $K=FD\oplus H(2;\un{1})^{(2)}$ where either
$D\in\{0,D_H(x_1^{p-1}x_2^{p-1})\}$ or $p=5$ and
$D=x_1^{p-1}\partial_2$, and there is an automorphism $\sigma$ of
$\Der H(2;\un{1})^{(2)}$ such that
$\sigma(\overline{T})=Fz_1\partial_1\oplus Fz_2\partial_2$ where
$z_i\in\{x_i,1+x_i\}$, $i=1,2$;

\smallskip

\item[(5)] $S\ot\OO(m;\un{1})\subset\overline{T}+K\subset\Der
(S\ot\OO(m;\un{1}))$ where $S$ is one of ${\mathfrak sl}(2)$,
$W(1;\un{1})$, $H(2;\un{1})^{(2)}$ and $m=1,2$. Let $\pi_2$ be the
canonical projection from $\Der (S\ot\OO(m;\un{1}))\cong\big((\Der
S)\ot\OO(m;\un{1})\big)\rtimes {\rm Id}_S\ot W(m;\un{1})$ onto
$W(m;\un{1})$. Then $\pi_2(K)\cong W(1;\un{1})$ if $m=1$ and
$\pi_2(K)\cong H(2;\un{1})^{(\epsilon)}$ if $m=2$, where
$\epsilon=1,2$;

\smallskip

\item[(6)] $S\ot\OO(1;\un{1})\subset K\subset
\widehat{S}\ot\OO(1;\un{1})$, where $S$ is one of ${\mathfrak
sl}(2)$, $W(1;\un{1})$, $H(2;\un{1})^{(2)}$ and either
$S=\widehat{S}$ or $S=H(2;\un{1})^{(2)}$ and
$\widehat{S}=H(2;\un{1})^{(1)}$. Moreover, $\overline{T}=F(h_0\ot
1)\oplus F({\rm Id}_{\widehat S}\ot(1+x_1)\partial_1)$ for some
toral element $h_0\in S$;

\smallskip

\item[(7)] $S\subset \overline{T}+K\subset S_p$ where $S$ is one
of the nonrestricted Cartan type Lie algebras $W(1;\un{2})$,
$H(2;\un{1},\Phi(\tau))^{(1)}$, $H(2;\un{1};\Delta)$ or else
$$H(2;(2,1))^{(2)}\subset \overline{T}+K\subset H(2;(2,1))_p;$$

\item[(8)] $K$ is either a classical Lie algebra of type ${\rm
A}_2,$ ${\rm B}_2$ or ${\rm G}_2$ or one of the restricted Cartan
type Lie algebras $W(2;\un{1})$, $S(3;\un{1})^{(1)},$
$H(4;\un{1})^{(1)},$ $K(3;\un{1})$.
\end{itemize}
\end{theo}
\pf 1) If $L(\alpha,\beta)$ is solvable, then we are in case (1).
So assume from now that $L(\alpha,\beta)$ is nonsolvable.

Recall from \cite{PS4} that $\rad_T L(\alpha,\beta)$ denotes the
maximal $T$-invariant solvable ideal of $L(\alpha,\beta)$, and
$L[\alpha,\beta]\,=\,L(\alpha,\beta)/\rad_T L(\alpha,\beta)$.
Since $\rad_T L(\alpha,\beta)\,=\, I(\alpha,\beta)\cap
L(\alpha,\beta)$, we have that
$$L[\alpha,\beta]\,\cong\,\psi(L(\alpha,\beta))=\,K\hookrightarrow \overline{T}+K.$$
As in \cite{PS4} we denote by
$\widetilde{S}=\oplus_{i=1}^r\,\widetilde{S}_i$ the sum of all
minimal $T$-invariant ideals of $K=L[\alpha,\beta]\ne (0)$. Since
the Lie algebra $\overline{T}+K$ is semisimple, it acts faithfully
on $\widetilde{S}$. We will identify $\overline{T}+K$ with a Lie
subalgebra of $\Der \widetilde{S}$. As shown in \cite[Sect. 4]{PS4},
we have that $r\in\{1,2\}$. Moreover, if $r=2$, then we are in
case~(3); see \cite[Theorem 4.1]{PS4}.

\smallskip

\noindent 2) From now on assume that
$\widetilde{S}=\widetilde{S}_1$ is the unique minimal
$\overline{T}$-invariant ideal of $K=L[\alpha,\beta]$. Recall from
\cite{PS4} that $TR(\widetilde{S})\le TR(L[\alpha,\beta])\le 2$.

Suppose $TR(\widetilde{S})=2$. If the Lie algebra $\widetilde{S}$
is restrictable, then \cite[Theorem 4.2]{PS4} shows that
$\widetilde{S}=L[\alpha,\beta]$. Therefore, if $\widetilde{S}$ ia
a classical Lie algebra or a restricted Lie algebra of Cartan
type, then we are in case~(8). As explained in \cite[Sect. 4]{PS4}
there is a natural restricted homomorphism
$\Psi_{\alpha,\beta}\colon\,T+L(\alpha,\beta)_p\rightarrow \Der
\widetilde{S}$ which maps $T+L(\alpha,\beta)\subset
T+L(\alpha,\beta)_p$ onto $\overline{T}+K$.

Suppose $\widetilde{S}\cong{\cal M}(1,1)$. Then
$\widetilde{S}=\Der \widetilde{S}$. Choose a two-dimensional
nonstandard torus $\overline{T}'$ in $\widetilde{S}$. There exists
a torus $T'$ in the restricted Lie algebra $T+L(\alpha,\beta)_p$
such that $\ker \alpha\cap\ker \beta\subset T'$ and
$\Psi_{\alpha,\beta}(T')=\overline{T}'$. By construction, $T'$ is
then a nonstandard torus of maximal dimension in $L_p$. Since no
such tori can exist by our general assumption, we derive that
$\widetilde{S}\not\cong{\cal M}(1,1)$.

If $\widetilde{S}$ is a nonrestricted Cartan type Lie algebra,
then \cite[Theorem 1.1]{PS3} yields that $\widetilde{S}$ is one of
$W(1,\un{2})$, $H(2;\un{1};\Phi(\tau))^{(1)}$,
$H(2;\un{1};\Delta)$, $H(2;(2,1))^{(2)}$. Applying \cite[Theorem
4.2]{PS4} shows that we are in case~(7).

\smallskip

\noindent 3) It remains to consider the situation where $r=1$ and
$TR(\widetilde{S})=1$, which is ruled by \cite[Theorem 4.4]{PS4}.
Due of Theorem 2.1, case (1) of that theorem is our case (2).
Assuming case (2) of Theorem 4.4 in \cite{PS4}, part (b) of the
proof in {\it loc. cit.} shows that there exists an automorphism
$\sigma$ of the Lie algebra $\Der H(2;\un{1})^{(2)}$ such that
$\sigma(K)\,=\,H(2;\un{1})^{(2)}\oplus FD$ and
$\sigma(\overline{T})\,=\,Fz_1\partial_1\oplus Fz_2\partial_2$,
where $D$ is as in case (4) of the present theorem and
$z_i\in\{x_i,1+x_i\}$, $i=1,2$.

Now assume we are in case (3) of \cite[Theorem 4.4]{PS4}. Then
$\widetilde{S}=S\ot{\mathcal O}(1;\un{1})$ where $S$ is one of
${\mathfrak sl}(2)$, $W(1;\un{1})$, $H(2;\un{1})^{(2)}$, and
$\Psi_{\alpha,\beta}(T)$ is spanned by $h_0\ot 1$ and ${\rm
Id}\ot(1+x_1)\partial_1$ for some nonzero toral element $h_0\in
S$. Moreover, $K\subset (\Der S)\ot\OO(1;\un{1})$. To show that
this
is our case (6) we can assume that $S=H(2;\un{1})^{(2)}$ (in
the other two cases $\Der S\,=\ad S$ and there is nothing to
prove).

If $K\not\subset H(2;\un{1})^{(2) }\ot\OO(1;\un{1})$, then there is
a root $\mu\in\Gamma(L,T)\cap(\F_p\alpha+\F_p\beta)$ such that
$K(\mu)\not\subset H(2;\un{1})^{ (2) } \ot\OO(1;\un{1})$.
Restricting the composite
$$\bar{\psi}\colon\,L(\alpha,\beta)\stackrel{\psi}{\longrightarrow}
K\hookrightarrow(\Der S)\ot\OO(1;\un{1})\twoheadrightarrow \Der S
$$ to the $1$-section $L(\mu)$ and using the above description of
$\Psi_{\alpha,\beta}(T)$ it is easy to observe that
$K(\mu)\cong\bar{\psi}(L(\mu))$ is sandwiched between $S$ and $\Der
S$. Consequently, $\mu$ is a Hamiltonian root and $K(\mu)\cong
L[\mu]$. Theorem~2.1 now yields $K(\mu)\cong
H(2;\un{1})^{(\epsilon)}$ where $\epsilon\in\{1,2\}$. If
$\epsilon=2,$ then $K(\mu)=K(\mu)^{(\infty)}\subset \widetilde{S}$
contrary to our choice of $\mu$. Hence $\epsilon=1$. Since
$K(\mu)\not\subset H(2;\un{1})^{(2)} \ot\OO(1;\un{1})$, it must be
that $\bar{\psi}(L(\mu))= H(2;\un{1})^{(1)}$ by Theorem 2.1. This
proves that $K\subset H(2;\un{1})^{(1)}\ot\OO(1;\un{1}),$ as
claimed.

In case~(4) of \cite[Theorem 4.4]{PS4} there exist
$\mu,\nu\in\Gamma(L,T)$ such that $L[\mu,\nu]\cong {\cal M}(1,1)$.
But then, as before, $L_p$ contains a nonstandard torus of maximal
dimension, violating our general assumption. Case~(5) of
\cite[Theorem 4.4]{PS4} is included into case (5) of the present
theorem.

Finally, consider  case~(6) of \cite[Theorem 4.4]{PS4}. Let $h_0$
be a nonzero toral element of $S=H(2;\un{1})^{(2)}$ and choose a
torus $T'$ in $T+L(\alpha,\beta)_p$ with $\dim T'=\dim T$ and
$h_0\ot 1\in \Psi_{\alpha,\beta}(T')$, Note that $\ker
\alpha\cap\ker\beta\subset T'$ and $T'$ is standard by our general
assumption. Choose $\mu\in\Gamma(K,\Psi_{\alpha,\beta}(T'))$ with
$\mu(h_0\ot 1)=0$ and regard it as a linear function on $T'$
vanishing on $\ker \alpha\cap\ker\beta$. Then
$${\mathfrak c}_S(h_0)\ot \OO(2;\un{1})\subset K(\mu)\subset
\big({\mathfrak c}_{\Der S}(h_0)\ot\OO(2;\un{1})\big)\rtimes{\rm
Id}\ot\pi_2(K).$$ Let $\varrho$ be the canonical projection
$\big({\mathfrak c}_{\Der S}(h_0)\ot\OO(2;\un{1})\big)\rtimes{\rm
Id}\ot\pi_2(K)\twoheadrightarrow \pi_2(K)$. Note that $\varrho$ is
$T'$-equivariant and maps the $1$-section $K(\mu)$ onto
$\pi_2(K)$. By \cite[Theorem 4.4(6)]{PS4}, $\pi_2(K)$ is
sandwiched between $H(2;1)^{(2)}$ and $H(2;\un{1}).$ This implies
that $\varrho(K(\mu))$ is semisimple. As $\varrho(K(\mu))$ is a
homomorphic image of the $1$-section $L(\mu)$, it must be that
$\varrho(K(\mu))\cong L[\mu]$. Theorem 2.1 now yields
$\pi_2(K)=\varrho(K(\mu))\cong H(2;\un{1})^{(\epsilon)}$ where
$\epsilon\in\{1,2\}$, completing the proof. \qed

We have now established (for $p>3$) a refined version of
\cite[Theorem II.2]{St92} which is the corrected version of
\cite[Theorem 6.3]{St89b} (or, equivalently, of \cite[Theorem
1.15]{BOSt}). In what follows we will need a description of
$\overline{T}$ in the respective cases.
\begin{prop}
With the assumptions and notations of Theorem 2.2, $\overline{T}$
has the following properties in the respective cases:

\begin{itemize}
\item[(1)] $\overline{T}=(0)$; \smallskip

\item[(2)] $\overline{T}\subset S$ and $\dim \overline{T}=1$;

\smallskip

\item[(3)]$ \overline{T}=Fh_1\oplus Fh_2$ where $h_i\in S_i$,
$\,i=1,2$;

\smallskip

\item[(4)] $\overline{T}$ is conjugate to $Fz_1\partial_1\oplus
Fz_2\partial_2$ where $z_i\in\{x_1,1+x_i\}$, $\,i=1,2$;

\smallskip

\item[(5)]$\overline{T}=F(h\ot 1)\oplus F(d\ot1+{\rm Id}\ot t)$
where $h\in S$ and $t\in \pi_2(K)$ are nonzero toral elements,
$d\in\Der S$ is toral, and $[d,h]=0$;

\smallskip

\item[(6)] $\overline{T}=F(h\ot 1)\oplus F{\rm
Id}\ot(1+x_1)\partial_1$ where $h\in S\setminus (0);$ and
$h^{[p]}=h$;

\smallskip

\item[(7)] $\overline{T}\subset S_p$ and, moreover,
$\overline{T}\cap S=(0)$ when $S\cong
H(2;\un{1};\Phi(\tau))^{(1)}$ and $\dim \overline{T}\cap S=1$
otherwise;

\smallskip

\item[(8)] $\overline{T}\subset K$.
\end{itemize}
\end{prop}
\pf (a) If $K=(0)$, then $T+L(\alpha,\beta)$ is solvable, hence
coincides with $I(\alpha,\beta)$. Therefore, $\overline{T}=(0)$.

\smallskip

\noindent (b) From now on suppose $K\ne (0)$. Then $(0)\ne
\overline{T}+K=(T+L(\alpha,\beta))/I(\alpha,\beta)$ is semisimple,
hence acts faithfully on its socle $\widetilde{S}$. Note that
$$\widetilde{S}\subset (\overline{T}+K)^{(1)}\subset K\subset\Der\widetilde{S}.$$
We regard $\widetilde{T}+K$ as a Lie subalgebra of $\Der
\widetilde{S}$. The restricted homomorphism
$\Psi_{\alpha,\beta}\colon\, T+L(\alpha,\beta)_p\rightarrow
\Der\widetilde{S}$ introduced in \cite[Sect.~4]{PS4} then maps $T$
onto $\overline{T}$. It also maps $L(\alpha,\beta)_p\subset L_p$
onto the $p$-envelope of $K$ in $\Der\widetilde{S}$. The latter
$p$-envelope will be denoted by $K_p$. It contains the
$p$-envelope $\widetilde{S}_p$ of $\widetilde{S}$ in $\Der
\widetilde{S}$.

Since $\widetilde{S}\subset
\overline{T}+\widetilde{S}_p\subset\Der \widetilde{S}$, the
restricted Lie algebra $\overline{T}+\widetilde{S}_p$, is
centerless. In conjunction with \cite[Theorem~1.2.6(3)]{St04} this
shows that if $\overline{T}$ is a torus of maximal dimension in
$\overline{T}+K_p$, then $\overline{T}\cap \widetilde{S}_p$ is a
torus of maximal dimension in $\widetilde{S}_p$. In view of
\cite[Theorem~1.2.8(3) and Theorem~1.3.11(3)]{St04} we also have
that
$$1\le TR(\overline{T}+K_p)\le TR(T+L(\alpha,\beta)_p)\le
TR(L_p(\alpha,\beta))\le 2.$$

\noindent (c) Let $\overline{T}+K$ be as in case (2) of Theorem
2.2. Then $\overline{T}+K=S$ where
$S=\widetilde{S}=\widetilde{S}_p$ is one of ${\mathfrak sl}(2)$,
$W(1;\un{1})$, $H(2;\un{1})^{(1)}$ or else $\overline{T}+K$ is
sandwiched between $S=H(2;\un{1})^{(2)}$ and $H(2;\un{1})^{(1)}$.
In any event, $\overline{T}+K$ is semisimple, restricted, and has
absolute toral rank $1$. It follows that $\dim \overline{T}=1$ and
$\overline{T}$ is a torus of maximal dimension in
$\overline{T}+K=\overline{T}+K_p$. But then $S\cap\overline{T}$ is
a torus of maximal dimension in $S=\widetilde{S}_p$, by our
concluding remark in part~(b). Hence $\overline{T}\subset S$, by
dimension reasons.

\smallskip

\noindent (d) In case~(4) of Theorem~2.2 we have $\dim
\overline{T}=2$. In cases (3) and (5)--(8), we have that either
$TR(\widetilde{S})=2$ or the Lie algebra $\widetilde{S}$ is not
simple. From this it follows that in the remaining cases of
Theorem~2.2 the torus $\overline{T}=\Psi_{\alpha,\beta}(T)$ acts
on $\widetilde{S}$ as a two-dimensional torus of derivations.
Indeed, otherwise there would exist a root $\gamma\in\Gamma(L,T)$
such that $\overline{T}+K=\overline{T}+K(\gamma)$ and this would
imply that $\widetilde{S}=K(\gamma)^{(\infty)}$ is simple; see
Theorem~2.1. By our concluding remark in part~(b),
$\overline{T}\cap\widetilde{S}$ is a torus of maximal dimension in
$\widetilde{S}_p$.

\smallskip

\noindent (e) Suppose we are in case (3) of Theorem~2.2. Then
$\widetilde{S}=S_1\oplus S_2$, $TR(S)=2$, and
$\widetilde{S}=\widetilde{S}_p$. Since $\dim (\overline{T}\cap
\widetilde{S}_p)=TR(\widetilde{S})$ due to our concluding remark
in part~(b), it must be that  $T\subset S_1\oplus S_2$.

\smallskip

\noindent (f) If $\overline{T}+K$ is as in case~(4), then
$\overline{T}$ already has the required form in view of
Theorem~2.2.

\smallskip

\noindent (g) Suppose we are in case (5) of Theorem~2.2. Because
$S\ot\OO(m;\un{1})_{(1)}$ is a $[p]$-nilpotent ideal of
$\widetilde{S}=S\ot\OO(m;\un{1})$, it follows that
$$\widetilde{S}_p\,=\,S\ot\OO(m;\un{1})+S_p\ot F\,=\,\widetilde{S}$$
and $TR(\widetilde{S})=TR(S)=1$. The discussion in part~(b) now
yields $\overline{T}\cap\widetilde{S}_p=F\tilde{h}$ for some nonzero
toral element $\tilde{h}\in \widetilde{S}$. According to
\cite[Theorem~2.6]{PS2}, there is an automorphism $\varphi$ of
$\Der\widetilde{S}$ such that $\varphi(\overline{T})\subset (\Der
S)\ot F \rtimes {\rm Id}\ot W(m;\un{1})$. Since $\varphi$ preserves
the socle of $\Der \widetilde{S}$, it must be that
$$\varphi(\tilde{h})\in\big(S\ot
{\mathcal O}(m;\un{1})\big)\cap\big((\Der S)\ot F\big)\,=\,S\ot F.$$
In other words, $\varphi(\tilde{h})=h\ot 1$ for some  nonzero toral
element $h\in S$.

Now let $\tilde{t}$ be any (nonzero) toral element such that
$\overline{T}=F\tilde{h}\oplus F\tilde{t}$. Write
$\varphi(\tilde{t})=d\ot 1+{\rm Id}\ot t$ with $d\in\Der S$ and
$t\in W(m;\un{1})$. It is straightforward to see that $d$ and $t$
are both toral, and $[d,h]=0$.

It remains to show that $t$ is a nonzero element of $\pi_2(K)$. If
$t=0$, then $\overline{T}$ lies in $(\Der S)\ot\OO(m;\un{1})$,
that is $\pi_2(\overline{T})=(0)$. Since $\overline{T}$ is a torus
of maximal dimension in $\overline{T}+K_p$, we apply
\cite[Theorem~1.2.8(4)]{St04} to get $TR(\pi_2(K))=0$. But
Theorem~2.2 says that in the present case $\pi_2(K)$ is a
semisimple restricted Lie algebra of absolute toral rank one.
Hence $Ft$ is a torus of maximal dimension in
$\pi_2(\overline{T}+K_p)\,=\,Ft+\pi_2(K)$. As $\pi_2(K)$ is a
restricted ideal of $\pi_2(\overline{T}+K_p)$, we conclude that
$t\in \pi_2(K)$; see \cite[Theorem~1.2.6(3)]{St04}.

\smallskip

\noindent (h) If $\overline{T}+K$ is as in case~(6), then
$\overline{T}$ already has the required form in view of
Theorem~2.2.

\smallskip

\noindent (i) Suppose  $\overline{T}+K$ is as in case~(7). If $S$
is one of $W(1;\un{2})$, $H(2;\un{1};\Phi(\tau))^{(1)}$,
$H(2;\un{1},\Delta)$, then Theorem~2.2 says that
$\overline{T}\subset S_p$. Suppose $S=W(1;\un{2})$. As
$W(1;\un{2})$ has codimension $1$ in $W(1;\un{2})_p$, it must be
that $\overline{T}\cap W(1;\un{2})\ne (0)$. Since all weight
spaces of the $\overline{T}$-module $W(1;\un{2})$ are
one-dimensional, by \cite[Theorem~V.2(2)]{St92}, we have that
$\dim(\overline{T}\cap S)=1$ in this case. If
$S=H(2;\un{1};\Phi(\tau))^{(1)}$, then
\cite[Theorem~VII.3(3)]{St92} yields $\overline{T}\cap S=(0).$
Suppose $S=H(2;\un{1};\Delta)$. Then $H(2;\un{1};\Delta)_p\,=\,
F(x_1\partial_1+x_2\partial_2)\oplus H(2;\un{1};\Delta)$. Hence
$S$ has codimension $1$ in $S_p$, implying $\overline{T}\cap S\ne
(0)$. Suppose $\overline{T}\subset S$. As explained in
\cite[Lemma~4.14(1)]{PS3}, for example, $T$ has at least $p^2-2$
weights on $S$. Since $\dim S=p^2$, we then have
$$S\,=\,\overline{T}\oplus\bigoplus_{\gamma\in\Gamma(S,T)}
S_\gamma, \quad |\Gamma(S,\overline{T})|=p^2-2,\, \quad \dim
S_\gamma=1\quad\ (\forall\,\gamma\in\Gamma(S,\overline{T})).$$ But
then the subset
$\overline{T}\cup\big(\bigcup_{\gamma\in\Gamma(S,\overline{T})}
S_{\gamma}^{[p]}\, \big)$ spans a diagonalizable Lie subalgebra of
$S_p$ containing $\overline{T}$. Since $\overline{T}$ is a torus of
maximal dimension in $S_p$, we get
$\bigcup_{\gamma\in\Gamma(S,\overline{T})} S_{\gamma}^{[p]}
\subset\overline{T}.$ However, this would imply that $S$ is a
restrictable Lie algebra, which is false. We conclude that $\dim
(\overline{T}\cap S)=1$ in the present case.

Now consider the case where $S=H(2;(2,1))^{(2)}$. It is well-known
(and easily seen) that $S_p=F\partial_{1}^p\oplus S$ and
$H(2;(2,1))_p=F\partial_{1}^p\oplus H(2;(2,1)).$ But then it is
clear that the restricted Lie algebra $H(2;(2,1))_p/S_p$ is
$[p]$-nilpotent. As $\overline{T}\subset H(2;(2,1))_p$ by
Theorem~2.2, we now derive that $\overline{T}\subset S_p$. Since
$S$ has codimension $1$ in $S_p$, we have $\overline{T}\cap S\ne
(0)$, while part~(f) of \cite[(10.1.1)]{BW} yields
$\overline{T}\not\subset S$. Therefore, $\dim(\overline{T}\cap
S)=1.$

\smallskip

\noindent (j) Finally, suppose $\overline{T}+K$ is as in case~(8)
of Theorem~2.2. Then $K=\widetilde{S}=\widetilde{S}_p$ is a
restricted simple Lie algebra of absolute toral rank $2$. Hence
$\dim(\overline{T}\cap\widetilde{S})=2$ by our discussion at the
end of part~(b). In other words, $\overline{T}\subset K$.\qed

\section{Optimal tori}
Given a torus $T$ of maximal dimension in $L_p$ we denote by
$\Gamma_p(L,T)$ the subset of all proper roots in $\Gamma(L,T)$.
Note that if $\gamma\in\Gamma(L,T)$ is proper, then
$\Gamma(L,T)\cap\F_p^*\gamma\subset\Gamma_p(L,T)$. We say that $T$
is an {\it optimal} torus if the number
$$r(T):=|\Gamma(L,T)\setminus \Gamma_p(L,T)|$$ is minimal possible; see
\cite[p.~242]{PS3}.

From now on we fix an $\F_p$-linear map $\xi\colon F\rightarrow F$
such that $\xi^p-\xi={\rm Id}_F$ (such a map is unique up to an
$\F_p$-valued additive function on $F$). The process of toral
switching (based on the ideas of \cite{Wi}, \cite{Wil83},
\cite{P86b}) has been described in detail in \cite[Sect.~2]{PS2}.
Given $x\in L_\alpha$, where $\alpha\in\Gamma(L,T)$, we denote by
$T_x$ the linear span of all
$t_x:=t-\alpha(t)(x+x^p+\cdots+x^{p^{m-1}})$ with $t\in T$ (here
$m=m(x)$ is the smallest nonnegative integer with the property that
$x^{p^{m}}\in T$). By \cite{P86b}, $T_x$ is a torus of maximal
dimension in $L_p$ and
$\Gamma(L,T_x)\,=\,\{\gamma_x\,|\,\,\gamma\in\Gamma(L,T)\}$, where
$$\gamma_x(t_x)=\gamma(t)-(\xi\circ\gamma)(x^{p^m})\alpha(t)\qquad\
(\forall\,t\in T).$$ We say that $T_x$ is obtained from $T$ by the
{\it elementary switching} corresponding to $x$. By \cite{P86b},
there exists an invertible linear operator $E_x=E_{x,\xi}\in{\rm
GL}(L_p)$ such that ${\mathfrak c}_L(T_x)=E_x({\mathfrak c}_L(T))$
and $L_{\gamma_x}=E_x(L_\gamma)$ for all $\gamma\in\Gamma(L,T)$.
Moreover, $E_x$ is a lopynomial in $\ad x$. The operator $E_x$ is
often referred to as a {\it generalized Winter exponential}.

It is immediate from the explicit form of generalized Winter
exponentials that $L(\alpha_x)=L(\alpha)$ and
$L(\alpha_x,\beta_x)=L(\alpha,\beta)$ for all roots
$\beta\in\Gamma(L,T)$; see \cite[pp.~218--222]{PS2} for more
detail. Recall that the torus $T_x$ is standard by our general
assumption. Solvable and classical roots are proper by definition.
If $\alpha$ is Witt or Hamiltonian, then one can always find an
element $x\in\bigcup_{i\ne 0}L_{i\alpha}$ such that the root
$\alpha_x\in\Gamma(L,T_x)$ is proper.

Our main goal in this section is to show that if $T$ is an optimal
torus, then {\it all\,} roots in $\Gamma(L,T)$ are proper, and
describe the action of optimal tori on $2$-sections.
\begin{lemm}
Let $K=L[\alpha,\beta]$ and $\overline{T}$ be as in Theorem~2.2,
and suppose that we are not in case~(5) of that theorem. Let $u\in
L_\alpha$, and assume that $\alpha\not\in\Gamma_p(L,T)$ and
$\alpha_u\in\Gamma_p(L,T_u)$. Then
$|\Gamma_p(L(\alpha,\beta),T_u)|>|\Gamma_p(L(\alpha,\beta),T)|$.
\end{lemm}
\pf We denote by $\bar{u}$ the image of $u$ in $K$.

\smallskip

\noindent (1) Since $L(\gamma)\cap I(\alpha,\beta)\subset
Q(\gamma)$ for all $\gamma\in\Gamma(L,T)$, it suffices to show
that under the above assumptions we have
$|\Gamma_p(K,\overline{T}_u)|>|\Gamma_p(K,\overline{T})|$. If
$K=(0)$, then the root $\alpha$ is solvable, hence proper. So this
case is impossible.

\smallskip

\noindent (2) Let $K$ be as in case~(2) of Theorem~2.2. Then
$K=\psi(L(\mu))$ for some $\mu\in\F_p\alpha+\F_p\beta$. Since all
roots in $(\F_p\alpha+\F_p\beta)\setminus\F_p\mu$ are solvable,
hence proper, it must be that
$K=\psi(L(\alpha))=\psi(L(\alpha_u))$, implying
$|\Gamma(K,\overline{T}_u)|=|\Gamma_p(K,\overline{T}_u)|$. Since
$\alpha$ is improper, the result follows.

\smallskip

\noindent (3) Let $K$ be as in case~(3) of Theorem~2.2. Then
$\overline{T}=(\overline{T}\cap S_1)\oplus (\overline{T}\cap S_2)$
by Proposition~2.3. Hence there exist $\mu,\mu'\in\Gamma(L,T)$
such that $K=K(\mu)+K(\mu')$ and $[K_{i\mu},K_{j\mu'}]=0$ for all
$i,j\in\F_p^*$. All roots in $(\F_p\alpha+\F_p\beta)\setminus
(\F_p\mu\cup\F_p\mu')$ are solvable. Hence $\alpha\in
\F_p\mu\cup\F_p\mu'$. No generality will be lost by assuming that
$\alpha=\mu'$. Since $K_{i\mu_u}=E_{\bar{u}}(K_{i\mu})$ for all
$i\in\F_p^*$, the preceding remark yields
$$K(\mu_u)\,=\,E_{\bar{u}}(K(\mu))\,=\,{\mathfrak
c}_{K(\mu)}(\overline{T}_u)\oplus\textstyle{\bigoplus}_{i\in\F_{p}^*}\,
K_{i\mu}$$ (as $E_{\bar{u}}$ is a polynomial in $\ad\bar{u}$, we
have that $E_{\bar{u}}(y)=y$ for all $y\in \bigcup_{i\ne
0}\,K_{i\mu}$). It follows that $\mu$ is proper if and only if
$\mu_u$ is. Since $\alpha=\mu'$ is improper and $\mu'_u$ is proper,
the result follows.

\smallskip

\noindent (4) Let $K$ be as in case~(4) of Theorem~2.2. Then
$K=FD\oplus H(2;\un{1})^{(2)}$ and there is an automorphism
$\sigma$ of the Lie algebra $\Der H(2;\un{1})^{(2)}$ such that
$\sigma(\overline{T})$ is one of
\begin{eqnarray*}
T_0&=&Fx_1\partial_1\oplus Fx_2\partial_2,\\
T_1&=&F(1+x_1)\partial_1\oplus Fx_2\partial_2,\\
T_2&=&F(1+x_1)\partial_1\oplus F(1+x_2)\partial_2.
\end{eqnarray*}
Replacing $K$ by its isomorphic copy $\sigma(K)\subset
H(2;\un{1})$ we may assume that $\sigma=\text{Id}$. Theorem~III.4
in \cite{St92} describes the $1$-sections of $\Der
H(2;\un{1})^{(2)}$ relative to $\sigma(\overline{T})$, hence the
$1$-sections of $K$ relative to $\overline{T}$ (the deliberations
in \cite[Sect.~III]{St92} only require that $p>2$). It is
immediate from the description in \cite[Sect.~III]{St92} that all
improper roots in $\Gamma(K,\overline{T})$ are Witt (no root in
$\Gamma(K,\overline{T})$ is Hamiltonian by dimension reasons).

If $\sigma(\overline{T})=T_0$, then \cite[Theorem~III.5]{St92}
shows that all roots in $\Gamma(K,\overline{T})$ are proper. As
$\alpha\not\in\Gamma_p(L,T)$, this case is impossible.

Suppose $\overline{T}=T_1$. Let $\mu$ be any $\overline{T}$-root
of $K$ such that $\mu(x_2\partial_2)\not\in\{0,1\}$. Then in the
notation of \cite[Proposition~III.3]{St92} we have
$b\not\in\{0,1\}$. Hence $b\ge 2$, in which case that proposition
yields $H(2;\un{1})_\mu\subset H(2;\un{1})_{(0)}$. Thus if
$\mu(x_2\partial_2)\ne 0$, then there is a unique $i_0\in\F_p^*$
with $H(2;\un{1})_{i_0\mu}\not\subset H(2;\un{1})_{(0)}$.
Moreover, {\it loc. cit.} shows that $H(2;\un{1})_{i_0\mu}\cap
H(2;\un{1})_{(0)}$ has codimension $1$ in $H(2;\un{1})_{i_0\mu}$.
>From this it is easy to deduce that $\overline{T}$ normalizes a
solvable subalgebra of codimension $\le 1$ in $K(\mu)$. As a
consequence, $\mu$ is a proper root in $\Gamma(K,\overline{T})$.

Applying these deliberations to our improper root $\alpha$ we obtain
$\alpha(x_2\partial_2)=0$. Then $H(2;\un{1})(\alpha)$ is spanned by
$\{k(1+x_1)^{k-1}x_2\partial_2-(1+x_1)^k\partial_1\,|\,\,k\in\F_p\}$.
Therefore, $K(\alpha)$ is isomorphic to the Witt algebra
$W(1;\un{1})$. Since $\alpha_u$ is a proper root in $\Gamma(L,T_u)$,
the torus $\overline{T}_u=\overline{T}_{\bar{u}}$ must normalize
$H(2;\un{1})_{(0)}$. But then $\overline{T}_u$ is conjugate to $T_0$
under the automorphism group of $\Der H(2;\un{1})^{(2)}$. Our
remarks earlier in the proof now show that all roots in
$\Gamma(K,\overline{T}_u)$ are proper.

Suppose $\overline{T}=T_2$. Let $\mu$ be any root in
$\Gamma(K,\overline{T})$ and assume by symmetry that
$\mu((1+x_1)\partial_1)\ne 0$. Set
$$
a\,:=\,\frac{\mu((1+x_2)\partial_2)}{\mu((1+x_1)\partial_1)},\qquad\,
b\,:=\,\frac{\mu((1+x_1)\partial_1)-
\mu((1+x_2)\partial_2)}{\mu((1+x_1)\partial_1)}.
$$ Clearly, $a,b\in\F_p$. If $b=0$, then it is easy to see that
$H(2;\un{1})(\mu)$ is solvable. If $b\ne 0$, then the $1$-section
$H(2;\un{1})(\mu)$ is spanned by all
$k(1+x_1)^{k-1}(1+x_2)^{ka+b}\partial_2-(ka+b)(1+x_1)^k(1+x_2)^{ka+b-1}\partial_1$
with $k\in\F_p$. It follows from this description that
$H(2;\un{1})(\mu)\cong W(1;\un{1})$ and $\overline{T}$ does not
normalize the unique subalgebra of codimension $1$ in
$H(2;\un{1})(\mu)$. Since $H(2;\un{1})(\mu)\subset
H(2;\un{1})^{(2)}$, this implies that in the present case any root
in $\Gamma(K,\overline{T})$ is either solvable or improper.
Moreover, the inequality $|\Gamma_p(K,\overline{T})|\le p-1$
holds.

Since $\alpha_u$ is a proper Witt root in
$\Gamma(K,\overline{T}_u)$, the above discussion shows that the
torus $\overline{T}_u$ is not conjugate to $T_2$. Hence
$\overline{T}_u$ is conjugate to either $T_0$ or $T_1$. But then
$|\Gamma_p(K,\overline{T}_u)|>p-1$ by our remarks earlier in the
proof. The result follows.

\smallskip

\noindent (5) Case~(5) of Theorem~2.2 cannot occur by our general
assumption.

\smallskip

\noindent (6) Let $K$ be as in case~(6) of Theorem~2.2. Since
$\widetilde{S}$ is a restricted ideal of $\Der \widetilde{S}$, we
have that $t_{\bar{u}}-t\in\widetilde{S}$ for all
$t\in\overline{T}$. Since $L(\alpha,\beta)$ is a $2$-section for
$T_u$, the pair $(K,\overline{T}_u)$ appears in Proposition~2.3(6).
Hence we may assume that
$$\overline{T}_u\,=\,F(h\ot 1)\,\textstyle{\oplus}\, F{\rm
Id}\ot(1+x_1)\partial_1.$$ For $k\in\F_p$ put $\widehat{S}_k:=\{x\in
\widehat{S}\,|\,\,[h,x]=kx\}$ and $S_k:=\widehat{S}_k\cap S$. Given
$\mu\in\Gamma(K,\overline{T}_u)$ set $a=a(\mu):=\mu(h\ot 1)$ and
$b=b(\mu):=\mu\big({\rm Id}\ot(1+x_1)\partial_1\big)$. Then
$K_\mu=K\cap\big(\widehat{S}_a\ot(1+x_1)^b\big)$.

If $a=a(\mu)=0$, then $K(\mu)\subset {\mathfrak
c}_{\widehat{S}}(h)\ot\OO(1;\un{1})$. As $Fh$ is a maximal torus of
$S$ and $\widehat{S}/S$ is solvable, $K(\mu)$ is solvable too. Then
$\mu$ is proper. If $a=a(\mu)\ne 0$, then $K(\mu)\,=\,
\bigoplus_{i=0}^{p-1}\,K\cap\big(\widehat{S}_{ia}
\ot(1+x_1)^{ib}\big)$.  The evaluation map ${\rm
ev}\colon\,\widehat{S}\ot\OO(1;\un{1})\twoheadrightarrow
\widehat{S}$ taking $x\ot f\in\widehat{S}\ot\OO(1;\un{1})$ to
$f(0)\cdot x\in\widehat{S}$ is a Lie algebra homomorphism.  It is
injective on $K(\mu)$ and the image ${\rm ev}(K(\mu))$ is sandwiched
between $S$ and $\widehat{S}$. From this it is immediate that when
$a=a(\mu)\ne 0$, the root $\mu$ is proper in
$\Gamma(K,\overline{T}_u)$ if and only if either $S={\mathfrak
sl}(2)$ or $S$ is one of $W(1;\un{1})$, $H(2;\un{1})^{(2)}$ and $h$
normalizes the standard maximal subalgebra of $S$.

Since $\alpha$ is improper and $u\in L_\alpha$, the $1$-section
$K(\alpha)=K(\alpha_u)$ is neither solvable nor classical. Hence
$a(\alpha_u)\ne 0$ and $S\ne {\frak sl}(2)$. The above discussion
now shows that $h$ normalizes the standard maximal subalgebra of
$S$. This, in turn, yields that all roots in
$\Gamma(K,\overline{T}_u)$ are proper. Consequently,
$|\Gamma(K,\overline{T}_u)|>|\Gamma(K,\overline{T})|$.

\smallskip

\noindent (7) Let $K$ be as in case~(7) of Theorem~2.2. If
$S=W(1;\un{2})$, then \cite[Sect.~V]{St92} shows that all roots in
$\Gamma(K,\overline{T}_u)$ are proper (it is proved in {\it loc.
cit.} that for any two-dimensional torus $\mathfrak t$ in $S_p$
either $\Gamma(S_p,{\mathfrak t})=\Gamma_p(S_p,{\mathfrak t})$ or
$\Gamma_p(S_p,{\mathfrak t})=\emptyset$). Thus the result holds in
this case.

If $S=H(2;\un{1};\Phi(\tau))^{(1)}$, then for any two-dimensional
torus $\mathfrak t$ in $S_p$ all roots in $\Gamma(S_p,{\mathfrak
t})$ are solvable, hence proper; see \cite[Theorem~VII.3]{St92}.
Since $\alpha$ is improper, this case cannot occur. If $S=H(2;{\un
1};\Delta)$, then $Fx_1\partial_1+Fx_2\partial_2$ is a toral
Cartan subalgebra of $S_p$. Applying \cite[Corollary~2.10]{PS2}
shows that $\overline{T}$ is a toral Cartan subalgebra of $S_p$ as
well. Now \cite[Proposition~4.9]{Wil83} gives the result. If
$S=H(2;(2,1))^{(2)}$, then \cite[Lemma~10.1.1]{BW} applies and
gives the result.

\smallskip

\noindent (8) Finally, suppose $K$ is as in case~(8) of
Theorem~2.2. If $K$ is classical, then all roots in
$\Gamma(K,\overline{T})$ are classical, hence proper. Thus this
case cannot occur. If $K=W(2;\un{1})$, then
$Fx_1\partial_1+Fx_2\partial_2$ is a toral Cartan subalgebra of
$K$. But then $\overline{T}$ is a toral Cartan subalgebra of $K$
too; see \cite[Corollary~2.10]{PS2}. As before,
\cite[Proposition~4.9]{Wil83} gives the result. If $K$ is one of
$S(3;\un{1})^{(1)}$, $H(4;\un{1})^{(1)}$, $K(3;\un{1})^{(1)}$,
then the deliberations of \cite[(5.8)]{BW} apply and yield the
result. \qed
\begin{lemm}
Let $K=L[\alpha,\beta]$ and $\overline{T}=F(h\ot 1)\oplus F(d\ot
1+{\mathrm Id}\ot t)$ be as in case~(5) of Theorem~2.2 and
Proposition~2.3.  If $S$ is of Cartan type, let $S_{(0)}$ denote
the standard maximal subalgebra of $S$. If $S={\frak sl}(2)$, put
$S_{(0)}=S$. Let $u\in L_\alpha$ and
$\mu\in\Gamma(K,\overline{T})$. Then the following hold:

\smallskip

\begin{itemize}
\item[(1)] If $\mu(h\ot 1)=0$, then $\mu$ is proper if and only if
$t\in W(m;\un{1})_{(0)}$.

\smallskip

\item[(2)]
 Suppose $\mu(h\ot 1)\ne 0$. If $t\not\in W(m;\un{1})_{(0)}$, then $\mu$
is proper if and only if $h\in S_{(0)}$. If  $t\in
W(m;\un{1})_{(0)}$, then $\mu$ is proper if and only if the torus
$\overline{T}$ normalizes the maximal compositionally classical
subalgebra of $S(\mu)\cong
\widetilde{S}(\mu)/\big(S\ot\OO(m;\un{1})_{(1)}\big)(\mu)$.

\smallskip

\item[(3)] If $h\in S_{(0)}$ and $t\in W(m;\un{1})_{(0)}$, then
all roots in $\Gamma(K,\overline{T})$ are proper.

\smallskip

\item[(4)] Assume that $\alpha\not\in\Gamma_p(L,T)$ and
$\alpha_u\in\Gamma_p(L,T_u)$. Then
$|\Gamma_p(L(\alpha,\beta),T_u)|>|\Gamma_p(L(\alpha,\beta),T)|$.
\end{itemize}
\end{lemm}
\pf In proving this lemma, it will be convenient to work with a
slightly more general two-dimensional torus $${\mathcal R}:=F(h\ot
1)\oplus F(D+{\mathrm Id}\ot t),$$ where $D\in (\Der
S)\ot\OO(m;\un{1})$, $h$ and $t$ are toral elements of $S$ and
$W(m;\un{1}),$ respectively, and $[D,h\ot 1]=0$. Recall that
$\widetilde{S}=S\ot\OO(m;\un{1})$ and $m\in\{1,2\}$.

\smallskip

\noindent (1) Let $\mu\in\Gamma(K,{\mathcal R})$ be such that
$\mu(h\ot 1)=0$. Then  $\widetilde{S}(\mu)\subset{\mathfrak
c}_S(h)\ot\OO(m;\un{1})$. Since $Fh$ is a maximal torus of $S$,
the subalgebra ${\mathfrak c}_S(h)$ is nilpotent. Consequently,
$\widetilde{S}(\mu)$ is a nilpotent ideal of $K(\mu)$. But then
the preimage of $\widetilde{S}(\mu)$ under the canonical
homomorphism $L(\alpha,\beta)\twoheadrightarrow K$ lies in the
radical of $L(\mu)$. Since  the subalgebra $\pi_2(K)^{(1)}$ of
$W(m;\un{1})$ is isomorphic to one of $W(1;\un{1})$,
$H(2;\un{1})^{(2)}$, it follows that $\mu\in\Gamma(L,T)$ is proper
if and only if $Ft$ normalizes the standard maximal subalgebra of
$\pi_2(K)^{(1)}$.

If $m=1$, then $\pi_2(K)=W(1;\un{1})$. In this case $\mu$ is proper
if and only if $t\in W(m;\un{1})_{(0)}$. Suppose $m=2$. Then the
subalgebra $\pi_2(K)^{(1)}\cap W(2;\un{1})_{(0)}$ has codimension
$\le 2$ in $\pi_2(K)^{(1)}$ and is either solvable or isomorphic to
${\mathfrak sl}(2)$ modulo its radical (because ${\mathfrak sl}(2)$
is the semisimple quotient of $W(2;\un{1})_{(0)}$). This shows that
$\pi_2(K)^{(1)}\cap W(2;\un{1})_{(0)}$ is the standard maximal
subalgebra of $\pi_2(K)^{(1)}\cong H(2;\un{1})^{(2)}$. Then again
$\mu\in\Gamma(L,T)$ is proper if and only if $t\in
W(m;\un{1})_{(0)}$.

\smallskip

\noindent (2) Now let $\mu\in\Gamma(K,{\mathcal R})$ be such that
$\mu(h\ot 1)\ne0$. Then $K_{i\mu}\subset\widetilde{S}$ for all
$i\in\F_p^*$, whence $K(\mu)={\mathfrak c}_K({\mathcal
R})+\widetilde{S}(\mu)$. As ${\mathfrak c}_K({\mathcal R})$ is
nilpotent, $K(\mu)^{(\infty)}=\widetilde{S}(\mu)^{(\infty)}$.
Combining this with Theorem~2.1(ii) we now derive that
$$L[\gamma]^{(\infty)}\cong\,
K(\gamma)^{(\infty)}/\big(K(\gamma)^{(\infty)}\cap\,\rad\,K(\gamma)\big)\,\cong\,
\widetilde{S}[\mu]^{(\infty)}.$$ This shows that
$\mu\in\Gamma_p(K,{\mathcal R})$ if and only if
$\mu\in\Gamma_p(\widetilde{S},{\mathcal R})$. Recall that the
evaluation map ${\mathrm
ev}\colon\,\widetilde{S}=S\ot\OO(m;\un{1})\twoheadrightarrow S$,
$\,\,\,x\ot f\mapsto f(0)\cdot x,$ is a Lie algebra homomorphism
whose kernel $S\ot \OO(m;\un{1})_{(1)}$ is a nilpotent ideal of
$\widetilde{S}$.

\smallskip

\noindent (a) Suppose $t\not\in W(m;\un{1})_{(0)}$. Then
\cite[Theorem~2.6]{PS2} shows that there exists $\sigma \in{\mathrm
Aut}\,\widetilde{S}$ such that $\sigma({\mathcal R})\,=\,F(h'\ot
1)\oplus F{\mathrm Id}\ot (1+x_1)\partial_1$ for some nonzero toral
$h'\in S$. Since ${\mathcal R}\cap\widetilde{S}=F(h\ot 1)$, we can
assume (after rescaling $h'$ if necessary) that $\sigma(h)=h'$. Set
$\mu':=\mu\circ\sigma^{-1}$, an element in
$\Gamma(\widetilde{S},\sigma({\mathcal R})).$ Since
$\widetilde{S}(\mu')\,=\,\sigma(\widetilde{S}(\mu))$, we have that
$\mu\in\Gamma_p(\widetilde{S},{\mathcal R})$ if and only if
$\mu'\in\Gamma_p(\widetilde{S},\sigma({\mathcal R}))$.

Put $a:=\mu'(h'\ot 1)$ and $b:=\mu'({\mathrm
Id}\ot(1+x_1)\partial_1)$. Note that $a\in\F_p^*$ by our present
assumption on $\mu$. Let $u\in S$ be such that $[h',u]=ru$. Then
$r\in\F_p$ (for $h'$ is toral) and $u\ot (1+x_1)^{rb/a}\in
\widetilde{S}_{(r/a)\mu'}$. From this it is immediate that the
evaluation map takes $\widetilde{S}(\mu')$ onto $S$. Since
$\widetilde{S}(\mu')\cap \ker {\mathrm
ev}\subset\rad\,\widetilde{S}(\mu')$, the maximal compositionally
classical subalgebra of $\widetilde{S}(\mu')$ is mapped onto that
of $S$. It follows that
$\mu'\in\Gamma_p(\widetilde{S},\sigma({\mathcal R}))$ if and only
if $h'\in S_{(0)}$. Since the subalgebra $S_{(0)}$ is invariant
under all automorphisms of $S$, we obtain that in the present case
$\mu\in\Gamma_p(L,T)$ if and only if $h\in S_{(0)}$.

\smallskip

\noindent (b) Now suppose $t\in W(m;\un{1})_{(0)}$. Then ${\mathcal
R}$ preserves the ideal $\ker\,{\mathrm ev}$. As a consequence, the
evaluation map induces a natural Lie algebra homomorphism
$\Phi\colon {\mathcal R}+\widetilde{S}\rightarrow \Der S$. Since
$\widetilde{S}(\mu)\cap \ker \Phi\subset\rad\,\widetilde{S}(\mu)$,
it is straightforward that $\mu\in\Gamma_p(\widetilde{S},{\mathcal
R})$ if and only if $\Phi({\mathcal R})$ normalizes the maximal
compositionally classical subalgebra of
$\Phi(\widetilde{S}(\mu))\subset S$. Since
$\Phi(\widetilde{S}(\mu))\cong
\widetilde{S}(\mu)/\big(S\ot\OO(m;\un{1})_{(1)}\big)(\mu)$, part~(2)
follows.

\smallskip

\noindent (3) Next assume that $h\in S_{(0)}$ and $t\in
W(m;\un{1})_{(0)}$, and suppose that $\mu$ is an  improper root in
$\Gamma(K,{\mathcal R})$. Part~(1) shows that $\mu(h\ot 1)\ne 0$,
while part~(2b) says that $\Phi(\mathcal R)$ does not normalize
the maximal compositionally classical subalgebra of $S(\mu)\cong
\widetilde{S}(\mu)/\big(S\ot\OO(m;\un{1})_{(1)}\big)(\mu)$. In
particular, $\widetilde{S}_{j\mu}\not\subset
S\ot\OO(m;\un{1})_{(1)}$ for some $j\in\F_p^*$. Note that
$h\in\Phi({\mathcal R})$. If $\Phi({\mathcal R})=Fh$, we have that
$S=S(\mu)$. But then $h\in S_{(0)}$ normalizes the maximal
compositionally classical subalgebra of $S(\mu)$, contrary to our
choice of $\mu$.

Therefore, $\dim \Phi({\mathcal R})=2$. This implies that
$S=H(2;\un{1})^{(2)}$. Since $\Phi({\mathcal R})\cap S=Fh\subset
S_{(0)}$, it follows from \cite[Proposition~III.1(5)]{St92} that
$\Phi({\mathcal R})$ is conjugate under ${\rm Aut}\,S$ to the
torus $Fx_1\partial_1\oplus Fx_2\partial_2$. But then all roots in
$\Gamma(S,\Phi({\mathcal R}))$ are proper by
\cite[Theorem~III.5]{St92}, contrary to part~(2b) and our choice
of $\mu$. This contradiction proves that all roots in
$\Gamma(K,{\mathcal R})$ are proper.

\smallskip

\noindent (4) We now apply the above results to $\overline{T}$.
Assume that $\alpha\not\in\Gamma_p(L(\alpha,\beta),T)$ and
$\alpha_u\in\Gamma_p(L(\alpha,\beta),T_u)$, where $u\in L_\alpha$.
Let $\bar{u}$ denote the image of $u$ in $K=L[\alpha,\beta]$.
Regard $\alpha$ as a $\overline{T}$-root of $K$. Our assumption on
$\alpha$ and $u$ implies that $\bar{u}\ne 0$ and
$\F_p^*\alpha\subset\Gamma(K,\overline{T})$.

\smallskip

\noindent (a) Suppose $\alpha$ vanishes on $\overline{T}\cap
\widetilde{S}$. As $\alpha$ is improper, part~(1) yields $t\not\in
W(m;\un{1})_{(0)}$. According to \cite[Theorem~2.6]{PS2}, we can
assume without loss of generality that $\overline{T}\,=\,F(h'\ot
1)\oplus F{\rm Id}\ot(1+x_1)\partial_1$. If $h'\not\in S_{(0)}$,
then parts~(1) and (2) yield $\Gamma_p(K,\overline{T})=\emptyset$.
Then $|\Gamma_p(K,\overline{T}_u)|\ge
p-1>|\Gamma_p(K,\overline{T})|$, forcing
$|\Gamma_p(L(\alpha,\beta),T_u)|>|\Gamma_p(L(\alpha,\beta),T)|$.

So from now we assume that $h'\in S_{(0)}$. Since $\alpha$
vanishes on $\overline{T}\cap \widetilde{S}$, we have
$\alpha(h'\ot 1)=0$. Note that the torus $\overline{T}_u$ is
spanned by the elements $(h'\ot 1)_{\bar{u}}$ and $\big({\rm
Id}\ot(1+x_1)\partial_1\big)_{\bar{u}}$. As $\alpha(h'\ot 1)=0$,
our discussion at the beginning of this section shows that $(h'\ot
1)_{\bar{u}}=h'\ot 1$ and $\alpha_u(h'\ot 1)=0$. Besides,
$\big({\rm Id}\ot(1+x_1)\partial_1\big)_{\bar{u}}=D+{\rm Id}\ot
t'$ for some $D\in (\Der S)\ot\OO(m;\un{1})$ and $t'\in
W(m;\un{1})$. Since $\alpha_u$ is proper and vanishes on
$\overline{T}_u\cap\widetilde{S}=F(h'\ot 1)$, part~(1) shows that
$t'\in W(m;\un{1})_{(0)}$. But then part~(3) implies that all
roots in $\Gamma(K,\overline{T}_u)$ are proper, yielding
$|\Gamma_p(L(\alpha,\beta),T_u)|>|\Gamma_p(L(\alpha,\beta),T)|$.

\smallskip

\noindent (b) Suppose $\alpha(h\ot 1)\ne 0$. Then $\bar{u}\in
[h\ot 1,K_\alpha]\subset\widetilde{S}_\alpha$. If $t\not\in
W(m;\un{1})_{(0)}$, then, as before, it can be assumed that
$\overline{T}\,=\, F(h'\ot 1)\oplus F{\rm Id}\ot(1+x_1)\partial_1$
and $\alpha(h'\ot 1)\ne 0$; see \cite[Theorem~2.6]{PS2}. Since
$\alpha$ is improper, part~(2) shows that $h\not\in S_{(0)}$. But
then all roots in $\Gamma(K,\overline{T})$ are improper, by (1)
and (2), implying
$|\Gamma_p(L(\alpha,\beta),T_u)|>|\Gamma_p(L(\alpha,\beta),T)|$.

So assume $t\in W(m;\un{1})_{(0)}$. Then all roots in
$\Gamma(K,\overline{T})$ vanishing on $h'\ot 1$ are proper by (1).
Let $\Phi\colon\,\overline{T}+\widetilde{S}\rightarrow \Der S$ be
the Lie algebra homomorphism from part~(2b), so that
$\Phi(\widetilde{S})=\big(S\ot\OO(m;\un{1})\big)/\big(S\ot\OO(m;\un{1})_{(1)}\big)
\cong S$, $\Phi(h'\ot 1)=h'$, and $\Phi(d\ot 1+{\rm Id}\ot t)=d$.
As $\ker\Phi$ is solvable, a root $\nu\in\Gamma(K,\overline{T})$
with $\nu(h'\ot 1)\ne 0$ is proper if and only if
$\Phi(\overline{T})$ normalizes the maximal compositionally
classical subalgebra of the $1$-section $S(\nu)\cong
\widetilde{S}(\nu)/\big(S\ot\OO(m;\un{1})_{(1)}\big)(\nu)$. Since
$\alpha$ is improper, the above discussion shows that
$\widetilde{S}(\alpha)\not\subset S\ot\OO(m;\un{1})_{(1)}$.

Since $\bar{u}\in\widetilde{S}_\alpha$ and $S$ is restrictable, we
have that $(h'\ot 1)_{\bar{u}}\in \widetilde{S}$ and  $\big({\rm
Id}\ot t\big)_{\bar{u}}=D'+{\rm Id}\ot t$ for some $D'\in (\Der
S)\ot\OO(m;\un{1})$. Set $h'':=\Phi((h'\ot 1)_{\bar{u}})$, the
image of $(h'\ot 1)_{\bar{u}}$ in
$\big(S\ot\OO(m;\un{1})\big)/\big(S\ot\OO(m;\un{1})_{(1)}\big)$.
It is immediate from the proof of \cite[Lemma~2.5]{PS2} that
$\overline{T}_u$ is conjugate under ${\rm Aut}\, \widetilde{S}$ to
the torus $F(h''\ot 1)\oplus F(D''+{\rm Id}\ot t)$ for some
$D''\in (\Der S)\ot\OO(m;\un{1})$.

Suppose $\Phi(\overline{T})=Fh'$. Then
$S=S(\alpha)=S(\alpha_{\bar{u}}).$ Since $\alpha_{\bar{u}}$ is a
proper root, it must be that $h''\in S_{(0)}$. Combining the
preceding remark with (3) we now obtain that all roots in
$\Gamma(K,\overline{T}_u)$ are proper. Then
$|\Gamma_p(L(\alpha,\beta),T_u)|>|\Gamma_p(L(\alpha,\beta),T)|$.

Now suppose $\dim \Phi(\overline{T})=2.$ Then
$S=H(2;\un{1})^{(2)}$ and $\Phi(\overline{T})$ is conjugate under
${\rm Aut}\, S$ to one of the tori $T_i$, $i=0,1,2,$ from part~(4)
of the proof of Lemma~3.1. By our earlier remarks, $\alpha$ is an
improper root in $\Gamma(S,\Phi(\overline{T}))$. Therefore,
$\Phi(\overline{T})$ is not conjugate to $T_0$ and
$\alpha\in\Gamma(S,\Phi(\overline{T}))$ is Witt; see
\cite[Sect.~III]{St92} for more detail.

If $\Phi(\overline{T})$ is conjugate to $T_1$, then part~(4) of the
proof of Lemma~3.1 shows that $\overline{T}_u$ normalizes $S_{(0)}$.
Since the latter is invariant under all automorphisms of $S$, it
follows that $h''\in S_{(0)}$. But then all roots in
$\Gamma(K,\overline{T}_u)$ are proper by (3). If
$\Phi(\overline{T})$ is conjugate to $T_2$, then part~(4) of the
proof of Lemma~3.1 implies that there is a $\kappa\in
\Gamma(K,\overline{T})$ such that
$\Gamma_p(K,\overline{T})\subset\F_p^*\kappa$. Since $\alpha_u$ is a
proper Witt root in $\Gamma(S,\Phi(\overline{T}_u))$, part~(4) of
the proof of Lemma~3.1 also shows that
$|\Gamma_p(S,\Phi(\overline{T}_u))|>p-1$. But then
$|\Gamma_p(K,\overline{T})|\le p-1<
|\Gamma_p(S,\Phi(\overline{T}_u))|\le |\Gamma_p(K,\overline{T}_u)|$,
proving that
$|\Gamma_p(L(\alpha,\beta),T_u)|>|\Gamma_p(L(\alpha,\beta),T)|$ in
all cases. \qed

\begin{theo}
If $T$ is an optimal torus in $L_p$, then all roots in $\Gamma(L,T)$
are proper.
\end{theo}
\pf The proof of \cite[Proposition~10.4.1]{BW} applies without
changes, since for that proof one only needs the conclusions of
Lemmas~3.1 and 3.2. \qed
\begin{cor} Let $T$ be an optimal torus in $L_p$.
With the notations of Theorem 2.2 we have the following description
of $\overline{T}$ in the respective cases of that theorem:

\begin{itemize}
\item[(1)] $\overline{T}=(0)$.

\smallskip

\item[(2)] $\overline{T}\subset S$ is conjugate under an
automorphism of $S$ to $Fx_1\partial_1$ if $S=W(1;\un{1})$ and to
$F(x_1\partial_1-x_2\partial_2)$ if $S=H(2;\un{1})^{(2)}$.

\smallskip

\item[(3)]$ \overline{T}=Fh_1\oplus Fh_2$ where $h_i\in S_i$.
Moreover, for $i=1,2$, the torus $Fh_i$ is conjugate under ${\rm
Aut}\,S_i$ to $Fx_1\partial_1$ if $S_i=W(1;\un{1})$ and to
$F(x_1\partial_1-x_2\partial_2)$ if $S_i=H(2;\un{1})^{(2)}$.
\smallskip

\item[(4)] $\overline{T}$ is conjugate under an automorphism of
$S$ to the torus $Fx_1\partial_1\oplus Fx_2\partial_2$.

\smallskip

\item[(5)] Let $\{e_0,h_0,f_0\}$ be a standard basis of
${\mathfrak sl}(2)$. For $s=1,2$, let $y_1,\ldots, y_s$ be the
generating set of $\OO(s;\un{1})$ contained in
$\OO(s;\un{1})_{(1)}$, and let $D_1,\ldots,D_s\in W(s;\un{1})$ be
such that $D_i(y_j)=\delta_{ij}$ for all $1\le i,j\le s$. Then
$\overline{T}$ is conjugate under ${\rm Aut}\,\widetilde{S}$ to
one of the following tori:
\begin{eqnarray*}
&&\quad\ \,\,{\mathrm span}\big\{h_0\ot 1,\,{\rm Id}\ot
x_1\partial_1\big\}\mbox{ if }\, m=1 \mbox{ and }\, S={\mathfrak
sl}(2),\\
&&\quad\ \,\,{\mathrm span}\big\{y_1D_1\ot 1,\, {\rm Id}\ot
x_1\partial_1\big\} \mbox{ if }\, m=1 \mbox{ and }\,
S=W(1;\un{1}),\\
&&\quad\ \,\, {\mathrm span}\big\{(y_1D_1-y_2D_2)\ot
1,\,r(y_1D_1+y_2D_2)\ot 1+{\rm Id}\ot x_1\partial_1\big\}\ {\mathit
with}\ \, r\in\F_p,
\end{eqnarray*}
if  $m=1$  and $S=H(2;\un{1})^{(2)}$,
\begin{eqnarray*}
&&{\mathrm span}\big\{h_0\ot 1,\,{\rm Id}\ot
(x_1\partial_1-x_2\partial_2)\big\}\mbox{ if }\, m=2 \mbox{ and
}\, S={\mathfrak
sl}(2),\\
&&{\mathrm span}\big\{y_1D_1\ot 1,\, {\rm Id}\ot
(x_1\partial_1-x_2\partial_2)\big\} \mbox{ if }\, m=2 \mbox{ and
}\,
S=W(1;\un{1}),\\
&& {\mathrm span}\big\{(y_1D_1-y_2D_2)\ot 1,\,r(y_1D_1+y_2D_2)\ot
1+{\rm Id}\ot (x_1\partial_1-x_2\partial_2)\},
\end{eqnarray*}
$r\in\F_p,$ if  $m=2$  and $S=H(2;\un{1})^{(2)}$.
\smallskip

\item[(6)] $\overline{T}=F(h\ot 1)\oplus F{\rm
Id}\ot(1+x_1)\partial_1$ where $h\in S\setminus (0) $ and
$h^{[p]}=h$. If $S$ is of Cartan type, then $h\in S_{(0)}$.
\smallskip

\item[(7)] $\overline{T}\subset S_p$ and $\dim \overline{T}\cap
S=q$ where $q=0$ if $S = H(2;\un{1};\Phi(\tau))^{(1)}$ and $q=1$
otherwise;

\smallskip

\item[(8)] $\overline{T}\subset S_{(0)}$, where $S_{(0)}=S$ if $S$
is classical and $S_{(0)}$ is the standard maximal subalgebra of
$S$ if $S$ is of Cartan type.
\end{itemize}
\end{cor}
\pf (1) is clear.

\smallskip

\noindent (2) According to Proposition~2.3, one has
$\overline{T}\subset S$ and $\dim \overline{T}=1$. Therefore,
$S\cong L[\gamma]^{(\infty)}$ for some $\gamma\in\Gamma(L,T)$.
Suppose $S$ is of Cartan type. By Theorem~3.3, $\gamma$ is a proper
$T$-root. It follows that $\overline{T}\subset S_{(0)}$. The
statement now follows from Demushkin's theorem; see
\cite[(7.5)]{St04}. Part~(3) is analogous to (2). Part~(4) follows
from \cite[Theorem~III.5(1)]{St92}.

\smallskip

\noindent (5)  By Proposition~2.3, $\overline{T}$ is conjugate under
${\rm Aut}\,\widetilde{S}$ to $F(h\ot 1)\oplus F(d\ot 1+{\rm Id}\ot
t)$, where $h\in S$, $t\in\pi_2(K)\subset W(m;\un{1})$ are nonzero
toral elements, $d$ is a toral element of $\Der S$, and $[h,d]=0$.
Since all roots in $\Gamma(K,\overline{T})$ are proper, by
Theorem~3.3, Lemma~3.2 implies that $t\in W(m;\un{1})_{(0)}$ and
$h\in S_{(0)}$. Put ${\mathcal R}:=Fh+Fd$, a torus in $\Der S$.

If $\dim {\mathcal R}=1$, then $d$ is a scalar multiple of $h$; so
we can assume that $d=0$. In this case there exist $\phi\in{\rm
Aut}\, \OO(m;\un{1})$ satisfying $\phi\circ
H(2;\un{1})^{(2)}\circ\phi^{-1}=H(2;\un{1})^{(2)}$ if $m=2$ and
$\sigma\in{\rm Aut}\,S$ such that $\sigma(h)$ and $\phi\circ
t\circ \phi^{-1}$ are nonzero multiples of $h_0$ and
$x_1\partial_1$, $y_1D_1$ and $x_1\partial_1$, $y_1D_1-y_2D_2$ and
$x_1\partial_1$, $h_0$ and $x_1\partial_1-x_2\partial_2$, $y_1D_1$
and $x_1\partial_1-x_2\partial_2$, $y_1D_1-y_2D_2$ and
$x_1\partial_1-x_2\partial_2$ in the six respective cases (when
$S=H(2;\un{1})^{(2)}$ and $m=2$, one should also keep in mind
Demushkin's theorem mentioned above). Clearly, the desired
normalization can be achieved with the help of
$\sigma\ot\phi\in{\rm Aut}\,\widetilde{S}$. Note that $r=0$ when
$\dim {\mathcal R}=1$.

Now consider the case where $\dim {\mathcal R}=2$. Since all
$T$-roots of $L$ are proper, so are all ${\mathcal R}$-roots of
$S=H(2;\un{1})$; see Lemma~3.2(2). If $m=1$ (resp., $m=2$), then
$t$ is a nonzero toral element in $W(1;\un{1})_{(0)}$ (resp., in
${H(2;\un{1})^{(2)}}_{(0)}$). As before, there exists $\phi\in
{\rm Aut}\,\OO(m;\un{1})$ satisfying $\phi\circ
H(2;\un{1})^{(2)}\circ\phi^{-1}=H(2;\un{1})^{(2)}$ if $m=2$ such
that $\phi\circ t\circ \phi^{-1}=x_1\partial_1$ (resp., $\phi\circ
t\circ \phi^{-1}=x_1\partial_1-x_2\partial_2$) when $m=1$ (resp.,
$m=2$). According to \cite[Theorem~III.5(1)]{St92}, there exists
$\sigma\in{\rm Aut}\,S$ such that
$$\sigma({\mathcal R})=F(y_1D_1)\oplus F(y_2D_2),\qquad
\sigma(h)\in\F_p^*(y_1D_1-y_2D_2).$$ Subtracting a multiple of $h$
from $d$ if necessary we may assume that
$\sigma(d)=r(y_1D_1+y_2D_2)$ for some $r\in F^*$. Since
$\overline{T}$ is two-dimensional and contains $(d\ot 1+{\rm
Id}\ot t)^{p}=d^{[p]}\ot 1+{\rm Id}\ot t$, it can only be that
$r^p=r$, that is $r\in\F_p$. As before, the desired normalization
can now be achieved with the help of $\sigma\ot\phi\in{\rm
Aut}\,\widetilde{S}$.

\smallskip

\noindent Part~(6) is analogous to (2). Part~(7) has already been
proved; see Proposition~2.3.

\smallskip

\noindent (8) Assume that $S$ is a restricted Lie algebra of
Cartan type. As all $\overline{T}$-roots of $S$ are proper by
Theorem~3.3, the statement follows from the discussion in
\cite[Sect.~IX]{St92}. \qed

\smallskip

We now have generalized the main results of \cite{St89b} to the
case where $p>3$. In particular, we have classified all
$T$-semisimple quotients of $2$-sections that occur in $L$; see
Theorem~2.2. In fact, our list is more precise that the list in
\cite[Theorem~II.2]{St92} which is the revised version of
\cite[Theorem~6.3]{St89b}. For $p>3$, all roots with respect to
optimal tori in $L_p$ are proper (Theorem~3.3). We recall that our
definition of ``properness'' differs from that introduced by Block
and Wilson.

\section{The subalgebra $Q(L,T)$}

Our next goal is to show that all deliberations of \cite{BOSt} are
valid for $p>3$. Theorem~1.15 of \cite{BOSt} can now replaced by the
stronger Theorem~2.2 of the present work, and Theorem~1.16 of
\cite{BOSt} can be substituted by the stronger Corollary~3.4.

Lemmas~2.1--2.4 of \cite{BOSt} generalize easily to the case where
$p>3$: the proofs of Lemmas~2.1 and 2.2 require no changes, while
Lemmas~2.3 and 2.4  are even easier to prove now, as we acquired
more information on $\overline{T}$.

Let us look at the proof of \cite[Lemma~2.5]{BOSt}. Suppose $S\cong
W(1;\un{2})$. Since all roots are proper, $\overline{T}$ cannot be
as in \cite[Theorems~V.2 and V.3]{St92}. So $\overline{T}$ is as in
\cite[Theorem~V.4]{St92}, hence the result follows. Suppose $S\cong
H(2;\un{1};\Phi(\tau))^{(1)}$. Then \cite[Theorem~VII.3]{St92}
yields that all $1$-sections of $S$ are solvable. This is the claim.
Suppose $S\cong H(2;(2,1))^{(2)}$. Then Corollary~VI.3 and
Theorem~VI.4 of \cite{St92} imply the statement of
\cite[Lemma~2.5]{BOSt} in this case. Note that, apart from
straightforward computations, only \cite[(10.1.1)]{BW} is used in
\cite[Sect.~VI]{St92}, and that holds for $p>3$. When $S\cong
H(2;\un{1};\Delta)$, the proof of Lemma~2.5 in \cite{BOSt} relies
only on elementary observations and \cite[(11.1.3)]{BW}, which holds
for $p>3$.

The proof of Lemma~2.6 in \cite{BOSt} requires no changes (for a
more elaborate computation see \cite[Sect.~IX]{St92}). Corollary~2.7
of \cite{BOSt} holds for $p>3$ too.

Thus all results of \cite[Sect.~2]{BOSt} hold for $p>3$. It is now a
matter of routine to check that all results of \cite[Sect.~3]{BOSt}
remain valid for $p>3$ as well. As a consequence, we obtain the
following:
\begin{theo}
Suppose all tori of maximal dimension in $L_p$ are standard and
let $T$ be an optimal torus in $L_p$. If $\alpha\in\Gamma(L,T)$ is
classical or solvable, define $Q_\alpha:=L_\alpha$. If
$\alpha\in\Gamma(L,T)$ is Witt or Hamilton, denote by
$L[\alpha]_{(0)}$ the standard maximal subalgebra of $L[\alpha]$
and define $$Q_\alpha:= \,
\psi_\alpha^{-1}\big(L[\alpha]_{(0)}\big)\cap L_\alpha,$$ where
$\psi_\alpha\colon\,L(\alpha)\twoheadrightarrow L[\alpha]$ stands
for the canonical homomorphism. Then
$$Q=Q(L,T):=\,H\oplus\bigoplus_{\gamma\in\Gamma(L,T)} Q_\gamma$$
is a $T$-invariant subalgebra of $L$.
\end{theo}
Note that $L=Q$ if and only if all roots in $\Gamma(L,T)$ are
solvable or classical. As another consequence, we obtain:
\begin{theo}
Let $T$ be an optimal torus in $L_p$. Given
$\alpha,\beta\in\Gamma(L,T)$ set $Q(\alpha,\beta):=Q\cap
L(\alpha,\beta)$ and let ${\mathcal J}(\alpha,\beta)$ denote the
maximal solvable ideal of $T+Q(\alpha,\beta)$. Let
$\varrho_{\alpha,\beta}\colon T+Q(\alpha,\beta)\twoheadrightarrow
\big(T+Q(\alpha,\beta)\big)/{\mathcal J}(\alpha,\beta)$ be the
canonical homomorphism, and set
$M:=\varrho_{\alpha,\beta}(Q(\alpha,\beta))$. Then one of the
following hold:
\begin{itemize}

\smallskip

\item[({\rm A})] $\ \ M=(0)$; \smallskip

\item[({\rm B})] $\ \ M\cong \,{\mathfrak sl}(2)$;

\smallskip

\item[({\rm C})]$\ \ M\cong \,{\mathfrak sl}(2)\oplus {\mathfrak
sl}(2)$;

\smallskip

\item[({\rm D})] $\ \ M\cong \,{\mathfrak sl}(2)\ot\OO(1;\un{1})$;

\smallskip

\item[({\rm E})]$\ \ M \cong \,H(2;\un{1};\Phi(\tau))^{(1)}$;

\smallskip

\item[({\rm F})] $\ \ M$ is classical of type ${\mathrm A}_2$,
${\mathrm B}_2$ or ${\mathrm G}_2$.
\end{itemize}
\end{theo}
The main result of \cite{PS4} tells us that if $L=Q$, then $L$ is
either classical or of Cartan type. So from now on we will assume
that $Q(L,T)\ne L$ for any optimal torus $T\subset L_p$.

Our next goal is to show that $T+Q$ is a maximal subalgebra of
$T+L\subset L_p$. To do this we will need a result on roots in
$\Gamma(L,T)$ sticking out of the subalgebra $Q$. The union of all
such roots is denoted by $\Phi_-$ in \cite{St93}. The set
$\Phi_-\cap \,\Gamma(L(\alpha,\beta),T)$ is described in
\cite[Theorem~1.9]{St93} for all types of $2$-sections
$K=L[\alpha,\beta]$ that occur in Theorem~2.2. It should be
stressed here that the proof of Theorem~1.9 in \cite{St93} only
relies on \cite[(10.1.1)]{BW}, \cite[Lemma~2.5]{BOSt} and
\cite[III, IV.4, IV.5, V.4, VII, IX]{St92}. Therefore, it holds
for $p>3$.
\begin{theo}
$T+Q$ is a maximal subalgebra of $T+L$.
\end{theo}
\pf If $T+Q=T+L$, then all roots in $\Gamma(L,T)$ are solvable or
classical. Since this case has been excluded, $T+Q$ is a proper
subalgebra of $T+L$. Let $G$ be a subalgebra of $T+L$ such that
$T+Q\subsetneq G$. Note that $G$ is $T$-invariant. There exists a
nonsolvable, nonclassical root $\alpha\in\Gamma(L,T)$ such that
$G\cap L(\alpha)\ne Q(\alpha)$. As $Q(\alpha)$ is a maximal
subalgebra of $L(\alpha)$, it must be that $L(\alpha)\subset G$.

\smallskip

\noindent (a) Pick $\beta\in\Gamma(L,T)\setminus\F_p\alpha$ and
consider the $2$-section $L(\alpha,\beta)$ and its $T$-semisimple
quotient $K=L[\alpha,\beta]$. We let $\psi\colon
L(\alpha,\beta)\twoheadrightarrow K$ denote the canonical
homomorphism, and put $\overline{G}:=\,\psi(G\cap
L(\alpha,\beta))$.

We will now go through the eight cases of Theorem~2.2. If $K$ is as
in case~(6) or as in case~(7) with $S$ being one of $W(1;\un{2})$ or
$H(2;(2,1))^{(2)}$, then $\Gamma(K,\overline{T})$ contains a root
$\delta$ which vanishes on $\psi(H)$. Indeed, in case~(6) this is
easily deduced from Corollary~3.4(6), so assume that we are in
case~(7). Combining \cite[Theorem~1.9]{St93} with Theorem~2.2(7) one
observes that $\psi(H)$ normalizes the standard maximal subalgebra
$S_{(0)}=\psi(Q(\alpha,\beta))\cap K^{(1)}$. Since the subalgebras
$W(1;\un{2})_{(0)}$ and $H(2;(2,1))_{(0)}$ are restricted,
Theorem~2.2(7) yields that $TR(\psi(H),K)=1$. Hence
$\overline{T}\cap\psi(H)_p$ is spanned by a single toral element. On
the other hand, it follows from \cite{St92} that
$|\Gamma(K,\overline{T})|= p^2-1$ in the present case. This shows
that there is a root $\delta\in \Gamma(K,\overline{T})$ which
vanishes on $\psi(H)$.

The $1$-section $L(\delta)$ is nilpotent, hence $K(\delta)\subset
\psi(Q(\alpha,\beta))\subset\overline{G}$. By \cite[(10.1.1),
(11.1.1)]{BW} and \cite[Theorem~IV.5(3)]{St92}, the root space
$K_\delta$ contains an element $x$ with $\alpha(x^{[p]})\ne 0$. The
adjoint endomorphism $\ad x$ acts invertibly on each subspace
$\bigoplus_{i\in\F_p}K_{s\alpha+i\delta}$ with $s\in\F_p^*$. As
$L(\alpha)\cup L(\delta)\subset G$, this yields
$K\subset\overline{G}$.

\smallskip

\noindent (b) Now suppose that $K$ has type different from the ones
considered in part~(a). As $\alpha$ is Witt or Hamiltonian, it must
be that $\Gamma(K/\overline{G},\overline{T})\subset
\Phi_-\setminus\F_p\alpha$. The description of $\Phi_-\cap
\Gamma(L(\alpha,\beta),T)$ in \cite[Theorem~1.9]{St93} now entails
that $\Gamma(K/\overline{G},\overline{T})$ has one of the following
types:
$$\emptyset,\ \,\{\mu\},\ \,\{\pm \mu\},\ \,\{\mu,\nu\},\ \,\{\mu,\pm\nu\},$$
where $\mu$ and $\nu$ are $\F_p$-independent roots in
$\Gamma(L(\alpha\,\beta),T)\setminus\F_p\alpha$. Consequently, the
$\alpha$-string through each of the displayed roots contains no more
than two roots from $\Gamma(K/\overline{G},\overline{T})$. As $p\ge
5$, there exists $j_0\in\F_p^*$ such that $[K_{j_0\alpha},K]\subset
\overline{G}$. On the other hand, the derived subalgebra of
$K(\alpha)/\rad\, K(\alpha)$ is isomorphic to either $W(1;\un{1})$
or $H(2;\un{1})^{(2)}$; see Theorem~2.1(ii). From this it is
immediate that $K(\alpha)=\ker\tau+\rad\, K(\alpha)$, where $\tau$
denotes the representation of $K(\alpha)$ in $K/\overline{G}$
induced by the adjoint action of $L$.

We now interpret this information globally. Set
$N:=L(\alpha)^{(\infty)}$, a subalgebra of $G$. Let $\beta$ be any
$T$-root of $L$ with $\beta\not\in\F_p\alpha$. If
$K=L[\alpha,\beta]$ is as in part~(a) of this proof, then
$K=\overline{G}$, forcing $L(\beta)\subset G$. Otherwise,
$[N,L(\beta)]\subset G$ by the discussion above. As a result,
$[N,T+L]\subset G$. Applying \cite[Lemma~4.1]{BOSt} (which holds in
all characteristics), we now deduce that $N$ acts nilpotently on
$L$. But then $N$ acts nilpotently on itself contrary to the fact
that $N=N^{(1)}\ne (0) $. This contradiction shows that $G=T+L$ and
hence that $T+Q$ is a maximal subalgebra of $T+L$.\qed
\begin{cor}
If $Q=Q(L,T)$ is solvable, then $L\cong W(1;\un{n})$ for some
$\un{n}$ and $Q\subset L_{(0)}\cong W(1;\un{n})_{(0)}$.
\end{cor}
\pf By Theorem~4.3, $Q$ is a maximal $T$-invariant subalgebra of
$L$. As $Q$ is solvable, \cite[Corollary~9.2.13]{St04} says that
$L$ is one of ${\mathfrak sl}(2)$, $W(1;\un{n})$,
$H(2;\un{n};\Phi)^{(2)}$ (up to isomorphism). Moreover, the proof
of Corollary~9.2.13 in \cite{St04} shows that either $L={\mathfrak
sl}(2)$ or $Q$ is a Cartan subalgebra of $L$  or $L=W(1;\un{n})$
and $Q\subset W(1;\un{n})_{(0)}$. Since $Q\ne L$, the first
possibility cannot occur. If $Q$ is a Cartan subalgebra of $L$,
then every $1$-section of $Q$ relative to $T$ is nilpotent. But
then it immediate from the definition of $Q$ and Theorem~2.1(ii)
that every $1$-section of $L$ relative to $T$ is solvable.
However, in this case $Q=L$ which is impossible as $L$ is simple.
Thus, $L= W(1;\un{n})$ and $Q\subset W(1;\un{n})_{(0)}$ as stated.
\qed
\section{The associated graded algebra}
In this section we will go through \cite{St93} in order to extend
the results there to our present situation. All references to
\cite{St93} will underlined; for example, \un{Lemma~2.3} will
indicate that we refer to Lemma~2.3 of \cite{St93}. We will adopt
the notation of \cite{St93}; in particular, $\Phi=\Phi(L,T)$ will
denote the set $\Gamma(L,T)\cup\{0\}$.

\un{Lemma~1.1} and \un{Corollary~1.2} are valid in all
characteristics; see \cite[\S1.2]{St04}. \un{Theorem~1.3} holds for
$p>3$;  see Theorem~2.1(ii). \un{Theorem~1.4} is our Theorem~2.2,
and \un{Theorem~1.5} is covered by the stronger Corollary~3.4 of the
present paper. \un{Theorem~1.6} is easily deduced from our
Theorem~4.3 and Corollary~4.4, while \un{Theorem~1.7} is our
Theorem~4.2. \un{Lemma~1.8} is often referred to as Schue's lemma;
it holds all characteristics. \un{Theorem~1.9} holds for $p>3$; see
our discussion before the proof of Theorem~4.3. \un{Corollary~1.10}
follows from \un{Theorem~1.9}, hence remains valid in our present
setting. The proofs of \un{Theorem~1.11} and \un{Theorem~1.12} only
need \un{Theorems~1.4 and 1.7}. Therefore, these results remain true
for $p>3$. \un{Lemma~1.13} does not require any restrictions on $p$;
see \cite[Proposition~1.3.7 and Corollary~1.3.8]{St04}.
\un{Lemma~1.14} (which is a reformulation of \cite[(3.1.2)]{BW})
does not need any restrictions on $p$ either. Summarizing, all
results of \cite[Sect.~1]{St93} are valid for $p>3$.

The proofs of \un{Lemmas~2.1 and 2.2} work for $p>3$. Note that
\un{Lemma~2.2(4)} holds if $\mu\in\Delta$. Indeed, if $\mu$ vanishes
on $H$, then \cite[Theorem~3.1]{PS4} implies that $L(\mu)$ acts
triangulably on $L$. However, one has to make some changes on
pp.~15, 16 of \cite{St93}.

\smallskip
\noindent \un{\bf Lemma~2.3} (new parts).

\smallskip

\begin{itemize}

\item[(3)] {\it If $TR\big(H,L(\alpha,\beta,\kappa)\big)\ge 2$,
then $\dim L(\alpha,\beta,\kappa)/Q(\alpha,\beta,\kappa)\le 4p$.}

\smallskip

\item[(4)] {\it Assume that $\kappa\in \Phi_{[-1]}$. If
$\sum_{i,j\in\F_p}\,\dim
L_{\kappa+i\alpha+j\beta}/Q_{\kappa+i\alpha+j\beta}\ge 5p,$ then}

\smallskip

\begin{enumerate}
\item[(a)] $\alpha,\beta\in\Delta$;

\smallskip

\item[(b)] $\kappa(Q_{i\alpha+j\beta})\ne 0$ for all $i,j\in\F_p$;

\smallskip

\item[(c)] $L_{\kappa+i\alpha+j\beta}\not\subset Q$ {\it for all}
$i,j\in\F_p$;

\smallskip

\item[(d)] $\sum_{i,j\in\F_p}\dim
L_{\kappa+i\alpha+j\beta}/Q_{\kappa+i\alpha+j\beta}=p^2.$
\end{enumerate}
\end{itemize}

\smallskip
\pf (3) If $L_{i\kappa+j\alpha+k\beta}\subset Q$ for all nonzero
$(i,j,k)\in\F_p^3$, then we are done. Thus, replacing $\kappa$ by a
suitable root in $\F_p\kappa+\F_p\alpha+\F_p\beta$, we may assume
that $L_\kappa\not\subset Q$. Since
$TR\big(H,L(\alpha,\beta,\kappa)\big)\ge 2$, we may assume, renaming
$\alpha$ if necessary, that $\alpha$ and $\kappa$ are are linearly
independent on $H$. We take $\beta\in\Delta$ if
$TR\big(H,L(\alpha,\beta,\kappa)\big)= 2$, and we take $\beta$ to be
independent of $\alpha$ and $\kappa$ as functions on $H$ if
$TR\big(H,L(\alpha,\beta,\kappa)\big)\ge 3$.

Given $\ell\in\F_p$ put $\rho_\ell:=\alpha+\ell\beta$. Then in both
cases $\kappa$ and $\rho_\ell$ are $\F_p$-independent as functions
on $H$. Since $L(\alpha,\beta,\kappa) \,=\,\sum_{\ell\,\in\,\F_p}
L(\kappa,\rho_\ell)+L(\kappa,\beta),$ we have that
\begin{eqnarray*}\dim
L(\alpha,\beta,\kappa)/Q(\alpha,\beta,\kappa)&\le& \sum_{\
\ell\,\in\,\F_p}\big(\dim
L(\kappa,\rho_\ell)/Q(\kappa,\rho_\ell)-\dim
L(\kappa)/Q(\kappa)\big)\\
&+&\dim L(\kappa,\beta)/Q(\kappa,\beta).
\end{eqnarray*}
As $L_\kappa\not\subset Q$, the root $\kappa$ is either Witt or
Hamiltonian. It follows that $\kappa$ vanishes on $H\cap\rad_T\,
L(\kappa,\rho_\ell)\subset H\cap\rad\,L(\kappa)$ for all
$\ell\in\F_p$. If $\rho_\ell$ does not vanish on $H\cap\rad_T\,
L(\kappa,\rho_\ell)$, then
$L_{i\kappa+j\rho_\ell}\subset\rad_T\,L(\kappa,\rho_\ell)$ for all
$i\in\F_p$ and $j\in\F_p^*$. Then
$L[\kappa,\rho_\ell]=L[\kappa,\rho_\ell](\kappa)$ is of type~(2);
see Theorem~2.2. If both $\kappa$ and $\rho_k$ vanish on
$H\cap\,\rad_T\, L(\kappa,\rho_\ell)$, then they are linearly
independent as elements in
$\Gamma(L[\kappa,\rho_\ell],\overline{T})$, forcing
$TR\big(\psi(H),L[\kappa,\gamma_\ell]\big)=2$. In this situation
Theorem~2.2 shows that  $L[\kappa,\rho_\ell]$ cannot be be of
type~(1), (2) or (6). If $L[\kappa,\rho_\ell]$ is of type~(7)  and
$S$ is one of $W(1;\un{2})$, $H(2;(2,1))^{(2)}$,
$H(2;\un{1};\Phi(\tau))^{(1)}$, then one of the roots in
$\Gamma(L[\kappa,\rho_\ell],\overline{T})$ vanishes on $\psi(H)$;
see \cite[(V.4), (V5), (VI.4), (VII.4)]{St92}. Since in the present
situation $\kappa$ and $\rho_\ell$ are linearly independent on
$\psi(H)$, this is not the case.

Thus, each $2$-section $L[\kappa,\rho_\ell]$ must be of type~(2),
(3), (4), (5), (8) or of type (7) with $S=H(2;\un{1};\Delta)$.
Theorem~1.9 now implies the following: If $\kappa$ is Witt, then
only $\kappa$ and at most three more roots stick out of
$Q(\kappa,\rho_\ell)$. Also, at most $p$ roots stick out of
$Q(\kappa,\beta)$. Then $\dim
L(\alpha,\beta,\kappa)/Q(\alpha,\beta,\kappa)\le 3p+p=4p$. If
$\kappa$ is Hamiltonian, then only $\pm\kappa$ and at most two more
roots stick out of $Q(\kappa,\rho_\ell)$, whereas the number of
roots sticking out $Q(\kappa,\beta)$ is at most $2p$. In this case
we have $\dim L(\alpha,\beta,\kappa)/Q(\alpha,\beta,\kappa)\le
2p+2p=4p$. The claim follows.

\smallskip

\noindent (4) Let $n$ be the number of $2$-sections
$L(\kappa,\mu)$ with $\mu\in \F_p\alpha+\F_p\beta$ such that
$L[\kappa,\mu]$ is either of type~(6) or of type~(7) with  $S$
being $W(1;\un{2})$ or $H(2;(2,1))^{(2)}$. Our present assumption
on $\kappa$ is slightly weaker than that in the original
\un{Lemma~2.3(4)}. Arguing as in the original proof, we obtain
that $5p\le (3+n)p+3(1-n)$. Then $(n-2)p\ge 3(n-1)$, forcing
$n\not\in\{0,1,2\}$. Hence $n\ge 3$, and we can proceed as in the
original proof; see \cite[p.~16]{St93}. \qed

\smallskip

\un{Lemma~2.4} is an immediate consequence
\cite[Proposition~2.2]{St91b} which, in turn, relies on
\cite[Theorem~IV.3(1)]{St92}, some standard considerations, and a
version of our Corollary~3.4 (proved in \cite{St92}). Since
\cite[Theorem~IV.3]{St92} is true in all characteristics,
\un{Lemma~2.4} holds for $p>3$.

\smallskip

\noindent \un{\bf Lemma~2.5} (new proof). {\it Suppose
$\beta\in\Phi$ and $\alpha\in \Phi\setminus\Delta$, and put
$J:=J(Q,T)$. Then $[J_\alpha, Q_\beta]$ acts triangulably on $L$}.

\pf First suppose that $\alpha=i\mu$ and $\beta=j\mu$ for some
$\mu\in \Phi$ and $i,j\in\F_p$. Then
$\F_p^*\mu\cap\Delta=\emptyset$. If $i=-j$, then \un{Lemma~2.2(4)}
yields the assertion. If $i\ne -j$, the assertion follows from
\un{Lemma~2.2(3)}. So from now on we may assume that $\alpha$ and
$\beta$ are $\F_p$-independent.

Suppose the assertion is not true. Then there is $\kappa\in\Phi_-$
with $L_\kappa\not\subset Q$ and $\kappa([J_\alpha,Q_\beta])\ne 0$.
We have that $L_{\kappa+\alpha+\beta}
=[[J_\alpha,Q_\beta],L_\kappa]\not\subset Q$. \un{Lemma~2.3(1)} now
shows that $TR\big(H,L(\alpha,\beta,\kappa)\big)\ge 2$. The new
version of \un{Lemma~2.3(3)} then yields
$$\sum_{i,j\in\F_p^*} \dim L_{\kappa+i\alpha+j\beta}/Q_{\kappa+i\alpha+j\beta}\le 4p, $$
implying that the adjoint $(T+Q(\alpha,\beta))$-module
$\bigoplus_{i,j\in\F_p^*}\, L_{\kappa+i\alpha+j\beta}$ has a
composition factor of dimension $<p^2$. We call it $V$ and denote
by $\rho$ the corresponding representation of $T+Q(\alpha,\beta)$
in ${\mathfrak gl}(V)$.

In view of \un{Lemma~2.4}, $J\cap Q(\alpha,\beta)$ is a solvable
ideal of $T+Q(\alpha,\beta)$. If $(J\cap
Q(\alpha,\beta))^{(1)}\subset\ker \rho$, we set $I:=J\cap
Q(\alpha,\beta)^{(1)}+\ker \rho.$ Otherwise, choose $k\ge 1$
maximal subject to the condition $(J\cap
Q(\alpha,\beta))^{(k)}\not\subset\ker \rho$, and set $I:=(J\cap
Q(\alpha,\beta))^{(k)}+\ker\rho.$ Since
$[J_\alpha,Q_\beta]\not\subset \ker \rho$, it follows that in
either case $I$ is an ideal of $Q(\alpha,\beta)$ satisfying
$I^{(1)}\subset\ker\rho$ and $I\not\subset\ker \rho$.

It is easy to see that $I$ acts nilpotently on
$Q(\alpha,\beta)/\ker \rho$. As $Q_\beta\not\subset \ker\rho$ and
$Q_{\alpha+\beta}\not\subset\ker\rho$, we have that
$\alpha(I_\gamma)=\beta(I_\gamma)=0$ for all $\gamma\in\Phi$. By
general representation theory, there is a linear function $\chi$
on $I$ such that $\rho(u)-\chi(u){\rm Id}_V$ is a nilpotent
operator for all $u\in I$. It is immediate from the preceding
remark that $\chi(u)=\kappa(u)$ for all $u\in
\bigcup_{\gamma\in\Phi}\, I_\gamma$.

Parts~(a) and (b) of the original proof of \un{Lemma~2.5} do not
use any restriction on $p$. So let us take a look at part~(c):  If
there exists a nonzero $\lambda\in \Phi$ with
$I_\lambda\not\subset\ker \rho$, then part~(b) of the original
proof shows that $I_\lambda+\ker \rho$ is an ideal of
$T+Q(\alpha,\beta)$ acting nonnilpotenty on $V$. The above
discussion then yields that $\kappa(I_\lambda)\ne 0$. Since
$[I,Q(\alpha,\beta)]\subset\ker\rho$ by part~(b), it follows that
$$Q(\alpha,\beta)\subset\{u\in T+Q(\alpha,\beta)\,|\,\,\chi([I,u])=0\}
\subsetneq T+Q(\alpha,\beta).$$ By general representation theory,
the $(T+Q(\alpha,\beta))$-module $V$ is then induced from its
$Q(\alpha,\beta)$-submodule $V_0$ of dimension $\le p^{-1}\dim
V\le 4$.

Let $x_\beta\in Q_\beta$ be such that
$\kappa([J_\alpha,x_\beta])\ne 0$. Note that $J\cap\,
Q(\alpha,\beta)+Fx_\beta$ is a solvable subalgebra of
$Q(\alpha,\beta)$. Since $\dim V_0<p$, it must act triangulably on
$V_0$. Combining this with the Engel--Jacobson theorem one now
observes easily that $(J\cap
Q(\alpha,\beta)+Fx_\beta)^{(1)}\subset J$ acts nilpotently on $V$.
As a consequence, $[J_\alpha,x_\beta]$ acts nilpotently on $V$,
implying $\kappa([J_\alpha,x_\beta])=0$. This contradiction shows
that $I\subset H+\ker\rho$. The rest of the original proof works
for $p>3$.\qed

\smallskip

The original proofs of \un{Theorem~2.6}, \un{Corollary~2.7} and
\un{Corollary~2.8} go through for $p>3$ (in fact, the proof of
\un{Corollary~2.8} can be streamlined by making use of
\un{Corollary~2.7(3)}). The original proof of \un{Theorem~2.9} only
requires one very minor adjustment in the middle of p.~21 in
\cite{St93}:

\smallskip

\noindent ``If $\dim W\ge 5p$, then assertion (4a) of the new
\un{Lemma~2.3} yields $\alpha\in\Delta$, a contradiction. Thus
$\dim W<5p\le p^2$.''

\smallskip

Thus, all results in \cite[Sect.~2]{St93} essentially remain true
for $p=5,7$. Next, inspection shows that apart from standard
results valid in all characteristics, the arguments in
\cite[Sect.~3]{St93} rely only on results of \cite{St92} and
\cite{BOSt} valid for $p>3$, on \un{Theorems~1.9 and 2.6}, and on
\un{Lemma~2.3(5)} which we did not change. Hence all arguments and
results in \cite[Sect.~3]{St93} remain valid for $p>3$

In order to show that the proofs in \cite[Sect.~4]{St93}
generalize to the case where $p>3$ we first recall that all
results of \cite{St91a} are valid for $p>3$; see
\cite[Sect.~5]{PS4}. It should also be noted that
\cite[(4.7)]{BOSt} holds for $p>3$; see Corollary~4.4. Since
\cite[(2.4)]{St90} and \cite[(4.5)]{St91b} can now be substituted
by \cite[Theorem~D]{PS4}, inspection shows that all results used
in \cite[Sect.~4]{St93} are valid for $p>3$. Thus, what remains to
be  revised  in \cite[Sect.~4]{St93} is the proof of Claim~4 in
\un{Theorem~4.6} (this proof uses the inequality $p-1>5$ which is
no longer available in our situation).

\smallskip

\noindent {\it New proof.} Let $\overline{H}$ denote the image of
$H$ in $G_{\{0\}}$. As the $T$-root spaces of $G_{\{-1\}}$ are
$1$-dimensional, by part~(3) of the proof, the subalgebra
$\overline{H}$ is spanned by $\bar{h}_i$, $i=1,2$, and $\dim
\overline{H}=2$. By parts~(2) and (3) of the proof, we have that
$\bar{h}_2\in G_{\{0\}}^{(1)}$ and ${\rm tr}\, (\ad\!_{G_{\{-1\}}}
\bar{h}_1)\ne 0$.  It follows that $\overline{H}\cap
G_{\{0\}}^{(1)}\subset F\bar{h}_2$. Let $b$ denote the invariant
symmetric bilinear form on $G_{\{0\}}$ given by
$$b(x,y):=\,{\rm tr}\big(\ad\!_{G_{\{-1\}}}\,
x\circ\ad\!_{G_{\{-1\}}}\, y\big)\qquad \big(\forall\,x,y\in
G_{\{0\}}\big).$$ The radical $G_{\{0\}}^\perp:=\,\{x\in
G_{\{0\}}\,|\,\,b(x,G_{\{0\}})= 0\}$ of $b$ is an ideal of
$G_{\{0\}}$. Suppose $G_{\{0\}}^\perp\cap
G_{\{0\}}^{(1)}\cap\overline{H}\ne(0).$ Then $\bar{h}_2\in
G_{\{0\}}^\perp$, hence $b(\bar{h}_2,\bar{h}_2)=0$. Recall from
part~(3) of the proof that
$$G_{\{-1\}}\,=\,G_{\{-1\},\tilde{\beta}}+G_{\{-1\},-\tilde{\beta}}+\sum_{-1\le
i\le 1}G_{\{-1\},-\tilde{\alpha}+i\tilde{\beta}}$$ and
$\tilde{\alpha}(h_2)=0$, $\tilde{\beta}(h_2)=-1$. Since $p\ge 5$
and $\dim G_{\{-1\},\gamma} = 1$ for all weights $\gamma$ of
$G_{\{-1\}}$, this can only happen if
$G_{\{-1\}}=\,G_{\{-1\},\tilde{\alpha}}$. But then in the notation
of \cite[Sect.~4]{St93} we have $G=G_{(-1)}=Fz_0\oplus G'_{(0)}$.
Recall that $G'_{(0)}$ is the normalizer of $Fu_{p-2}\oplus J$,
where $J=\sum_{i\ge p-1}\,G_{[i]}$, and we are assuming that $J$
is a maximal ideal of $G$. From this it is immediate that the
subalgebra $G'_{(0)}/J$ of the simple Lie algebra $G/J$ acts
triangulably on $G/J$. However, $G'_{(0)}/J$ contains a copy of
${\mathfrak sl}(2)$ spanned by the images of $D(x^2)$, $D(xy)$ and
$D(y^2)$ in $G'_{(0)}/J$.

Therefore, $G_{\{0\}}^\perp\cap
G_{\{0\}}^{(1)}\cap\overline{H}=(0),$ implying
$G_{\{0\}}^\perp\cap G_{\{0\}}^{(1)}\,=\,\bigoplus_{\gamma\ne
0}\big(G_{\{0\}}^\perp\cap G_{\{0\}}^{(1)}\big)_\gamma$. Since all
elements in $\bigcup_{\gamma\ne 0} G_{\{0\},\gamma}$ act
nilpotently on $G_{\{-1\}}$, the Engel--Jacobson theorem yields
that $G_{\{0\}}^\perp\cap G_{\{0\}}^{(1)}$ acts nilpotently on
$G_{\{-1\}}$. Since $G_{\{-1\}}$ is an irreducible
$G_{\{0\}}$-module, $G_{\{0\}}^\perp\cap G_{\{0\}}^{(1)}=(0)$
necessarily holds. But then $b$ is nondegenerate on
$G_{\{0\}}^{(1)}$, forcing $G_{\{0\}}\,=\,G_{\{0\}}^{(1)}\oplus C$
where $C$ is the orthogonal complement to $G_{\{0\}}^{(1)}$ in
$G_{\{0\}}$. As $[G_{\{0\}}^{(1)},C]\subset C$ and $C^{(1)}\subset
G_{\{0\}}^{(1)}$, it must be that $C$ is a central ideal of
$G_{\{0\}}$. On the other hand, the center of $G_{\{0\}}$ has
dimension $\le 1$ by Schur's lemma, and $\bar{h}_1\not\in
G_{\{0\}}^{(1)}$ by our remarks earlier in the proof. Hence $C$
coincides with the center of $G_{\{0\}}$, and Claim~4 follows.
\qed

\medskip

We have already mentioned that Proposition~2.2 of \cite{St91b}
remains true for $p>3$. As a consequence, Lemma~2.4(1) of
\cite{St91b} is valid for $p>3$. Then one can see by inspection
that \un{Lemma~5.1} remains true for $p>3$ as well.

Part~(a) of the original proof of \un{Lemma~5.2} has to be
modified, however: in the course of the proof one has to show that
a certain torus $R$ is optimal, but the argument used in
\cite{St93} does not extend to the case where $p=5,7$. The
argument below will justify that $R$ is optimal for $p>3$.

We begin as in \cite[p.~46]{St93} and establish the existence of a
root $\kappa$ with
$$\kappa\not\in\Delta,\ \ \kappa(L_\alpha)\ne 0,\ \
\alpha([L_\kappa,L_{\beta-\kappa}]\ne 0.$$ Put
$K:=\,\bigcap_{\,n\ge 0}\,L(\alpha,\beta,\kappa)^{(n)}$ and let
$I$ be a maximal ideal of $K$. Put $G:=K/I$ and let $G_p$ be the
$p$-envelope of the simple Lie algebra $K$ in $\Der K$. Let
$\varphi\colon K\rightarrow G$ denote the canonical homomorphism.

\smallskip

\noindent {\it New part of the proof}. (a) Suppose $TR(G)=2$ and
put $N:=\varphi(K(\beta))$. Then one shows as in the original
proof that $N$ is a triangulable Cartan subalgebra of $G$ with
$TR(N,G)=TR(G)=2$. Let $N_p$ be the $p$-envelope of $N$ in $G_p$,
and let $R$ denote the unique maximal torus of $G_p$ contained in
$N_p$. Note that $\dim R=2$.

Suppose $L(\alpha,\beta)$ is not solvable. Then one shows as in
the original proof that $$L[\alpha,\beta]^{(1)}\cong
H(2;\un{1};\Phi(\tau))^{(1)}\cong
\big(L(\alpha,\beta)/\rad\,L(\alpha,\beta) \big)^{(1)}.$$ Our
choice of $\kappa$ then implies that $TR(K)=3$. Hence $I$ is a
$T$-invariant ideal of $K$. If $\beta(I_\alpha)\ne 0$ or
$\alpha(I_\beta)\ne 0$, then $K_{i\alpha+j\beta}\subset I$ for all
nonzero $(i,j)\in\F_p^2$ (one should keep in mind here that
$\beta(K_\alpha)\ne 0$ and $\alpha(K_\beta)\ne 0$). This yields
$TR(I)\ge 2$ forcing $TR(G)\le 1,$ a contradiction. Thus, it must
be that $\beta(I_\alpha)=\alpha(I_\beta)=0$. Therefore,
$$\beta(G_\alpha)\ne 0, \quad \alpha(G_\beta)\ne 0,\quad I\cap
K(\alpha,\beta)\subset H+\rad\,K(\alpha,\beta).$$ Since
$\alpha,\beta\in\Delta$ and $TR(G)=2$, this entails
$G=G(\alpha,\beta)$. Since $\kappa\not\in\Delta$, we now deduce that
$K_\kappa\subset I$. But then
$\alpha([I_\kappa,L_{\beta-\kappa}])=\alpha([K_\kappa,L_{\beta-\kappa}])=
\alpha([L_\kappa,L_{\beta-\kappa}])\ne 0$, showing that
$G_{i\alpha+j\beta}\subset I$ for all $i\in\F_p^*$ and $j\in\F_p$.
As this contradicts  our assumption that $TR(G)=2$, we deduce that
$L(\alpha,\beta)$ is solvable.

We now intend to show that all $R$-roots of $G$ are proper. The
$1$-sections of $G$ relative to $R$ are related to the
$2$-sections of $K$ relative  to $T$ as follows: Let $\mu$ be any
$T$-root of $K$. Since $I$ is $K(\beta)$-stable and $\dim\,R=2$,
the map $\varphi$ takes the subspace
$K_{\widetilde{\mu}}:=\,\bigoplus_{i\in\F_p}\,K_{\mu+i\beta}$ onto
a root space relative to $R$. Conversely, every $R$-root space of
$G$ is of the form $\varphi(K_{\widetilde{\mu}})$ for some
$\mu\in\langle\alpha,\beta,\kappa\rangle$. The nonzero roots
$\widetilde{\mu}\in\Phi(G,R)$ correspond to those
$\mu\in\Phi(K,T)$ which are $\F_p$-independent of $\beta$.

Let $\widetilde{\mu}$ be a nonzero root in $\Phi(G,R)$.  As $G$ is a
simple Lie algebra and $R$ is a torus of maximal dimension in $G_p$,
Theorem 2.1 shows that the derived subalgebra
$U:=G[\widetilde{\mu}]^{(1)}$ of
$G[\widetilde{\mu}]=G(\widetilde{\mu})/\rad\,G(\widetilde{\mu})$ is
either $(0)$ of one of ${\mathfrak sl}(2)$, $W(1;\un{1})$,
$H(2;\un{1})^{(2)}$. Furthermore, $U$ has codimension $\le 1$ in
$G[\widetilde{\mu}]$. If $U$ is either $(0)$ or ${\mathfrak sl}(2)$,
then $\widetilde{\mu}$ is solvable or classical, hence proper. So
from now on we may assume that $U$ is either $W(1;\un{1})$ or
$H(2;\un{1})^{(2)}$.

By analyzing the list of semisimple quotients in Theorem~2.2 one
finds out that $L[\mu,\beta]$ can only be of type $(2)$, $(3)$,
$(4)$, $(6)$ or $(7)$. Indeed, in case~(1) the Lie algebra
$L(\mu,\beta)$ is solvable, hence $U=(0)$. In case~(5) no
$1$-section in $L[\mu,\beta]$ is nilpotent, for otherwise one of the
roots in $\Gamma(L[\mu,\beta],\overline{T})$ would vanish on $h\ot
1\in\psi(H)$ (the notation of Proposition~2.3(6)). But then the
centralizer of $h\ot 1$ in $L[\mu,\beta]$ would be nilpotent,
contrary to the description in Theorem~2.2(6). In case~(8) the
inclusion $\overline{T}\subset \psi(H)$ holds as $L[\mu,\beta]$ is
simple and restricted. Since $U$ is of Cartan type, the Lie algebra
$L[\mu,\beta]$ cannot be classical. As a consequence, in both cases
(5) and (8) the set $\Phi(L[\mu,\beta],T)$ contains a nonzero
multiple of $\beta$. This, however, contradicts our assumption that
$\beta$ vanishes on $H$.

It is immediate from the definition of $K$ that $K(\mu,\beta)$ is an
ideal of $L(\mu,\beta)$ containing $\bigcap_{n\ge
0}\,L(\mu,\beta)^{(n)}$. Let $\pi\colon K(\mu,\beta)\rightarrow
L[\mu,\beta]$ denote the restriction to $K(\mu,\beta)\subset
L(\mu,\beta) $ of the canonical homomorphism $\psi\colon
L(\mu,\beta)\twoheadrightarrow L[\mu,\beta]$, and
$\varphi_{\mu,\beta}\colon K(\mu,\beta)\rightarrow G$ the
restriction to $K(\mu,\beta)\subset K$ of the epimorphism $\varphi$.
As explained above, $\varphi_{\mu,\beta}$ takes $K(\mu,\beta)$ onto
the $1$-section $G(\widetilde{\mu})$ with respect to $R$. Composing
$\varphi_{\mu,\beta}$ with the canonical homomorphism
$G(\widetilde{\mu})\rightarrow G[\widetilde{\mu}]$ we obtain a
surjective Lie algebra map $\nu\colon K(\mu,\beta)\twoheadrightarrow
G[\widetilde{\mu}]$.

Let $\widetilde{S}$ denote the socle of $L[\mu,\beta]$ and
$Q[\mu,\beta]$ the maximal compositionally classical subalgebra of
$L[\mu,\beta]$. Put
$$e:=\,\dim L[\mu,\beta]/Q[\mu,\beta].$$

Consider cases~$(2)$, $(4)$ and $(7)$. In each of these cases
$\widetilde{S}$ is a simple Lie algebra and the quotient
$L[\mu,\beta]/\widetilde{S}$ is nilpotent. Then
$\pi^{-1}(\widetilde{S})$ is simple modulo its radical and
$\pi(K(\mu,\beta))\supset \widetilde{S}$, by our earlier remarks.
Since $G[\widetilde{\mu}]$ is semisimple with simple socle
$U=G[\widetilde{\mu}]^{(1)}$, it must be that
$(\nu\circ\pi^{-1})(\widetilde{S})$ is either $(0)$ or $U$. As
$L[\mu,\beta]/\widetilde{S}$ is nilpotent, the first possibility
cannot occur. Therefore,
$$(\nu\circ\pi^{-1})(\widetilde{S})=\,U\cong\widetilde{S},\qquad\, \ker\nu\subset\ker\pi.$$
Consequently, $\ker\nu=\rad\,K(\mu,\beta)$. By our earlier
remarks, $\widetilde{S}$ is either $W(1;\un{1})$ or
$H(2;\un{1})^{(2)}$. We denote by $\widetilde{S}_{(0)}$ the
standard maximal subalgebra of $\widetilde{S}$.

As $TR(U)=1$, case~(7) is impossible.  It is easily seen that in
case~(2) we have $e=1$ if $\widetilde{S}$ is Witt and $e=2$ if
$\widetilde{S}$ is Hamiltonian. Lemma~2.2 of \cite{BOSt} holds for
$p>3$ and shows that in case~(4) we have $e=2$ and
$\widetilde{S}=H(2;\un{1})^{(2)}$. Furthermore,
$Q[\mu,\beta]\cap\widetilde{S}=\widetilde{S}_{(0)}$ in both cases.
Let $M:=(\nu\circ \pi^{-1})\big(Q[\mu,\beta])\big)$, a subalgebra
of $G[\widetilde{\mu}]$. Since $\ker \nu$ is solvable, $M$ has the
following properties:
\begin{itemize}
\item[1.]\, $M/\rad\,M$ is either $(0)$ or ${\mathfrak sl}(2)$;

\smallskip

\item[2.]\, $\dim\,G[\widetilde{\mu}]/M\le e$;

\smallskip

\item[3.]\, $M\cap U\subset U_{(0)}.$
\end{itemize}
Then $M$ coincides the maximal compositionally classical
subalgebra of $G[\widetilde{\mu}]$. Moreover, since $\beta$ is
solvable, it must be that $K(\beta)\subset
\pi^{-1}\big(Q[\mu,\beta]\big)$. This implies that $R$ normalizes
$M$. Consequently, $\widetilde{\mu}$ is a proper $R$-root.

Now suppose we are in case~(6). Then
$\widetilde{S}=S\ot\OO(1;\un{1})$, where $S$ is one of ${\mathfrak
sl}(2)$, $W(1;\un{1})$, $H(2;\un{1})^{(2)}$, and
$L[\mu,\beta]/\widetilde{S}$ is nilpotent. As $\widetilde{S}$ is
perfect and $\widetilde{S}/\rad\,\widetilde{S}$ is simple, this
gives $(\nu\circ\pi^{-1})(\widetilde{S})=S$ and $\ker\nu\subset
\ker\pi$. Hence $\ker\nu=\rad\,K(\mu,\beta)$.

We now proceed as before. Let $\mathfrak n$ denote the centralizer
of $S\ot\OO(1;\un{1})_{(p-1)}$ in $L[\mu,\beta]$. Since
$S\ot\OO(1;\un{1})\subset
L[\mu,\beta]\subset\widehat{S}\ot\OO(1;\un{1})$ by Theorem~2.2(6),
it is straightforward that ${\mathfrak n}=\rad\,L[\mu,\beta]$.
Then $\pi^{-1}({\mathfrak n})\subset \ker \nu$ by the preceding
remark. Using \cite[Lemma~2.4]{BOSt} (which holds for $p>3$) one
observes that when $S$ is Witt (resp., Hamiltonian), the
subalgebra $Q[\mu,\beta]+{\mathfrak n}$  has codimension $1$
(resp., $2$) in $L[\mu,\beta]$.  As $\pi^{-1}({\mathfrak
n})\subset\ker\nu$, it follows that
$(\nu\circ\pi^{-1})\big(Q[\mu,\beta]\big)$ is a maximal
compositionally classical subalgebra of $G[\widetilde{\mu}]$. It
contains $\nu(K(\beta))$ because $\beta$ is solvable. Therefore,
$\widetilde{\mu}$ is a proper $R$-root.

Finally, suppose we are in case~(3). Then
$L[\mu,\beta]/\widetilde{S}$ is nilpotent and
$\widetilde{S}=S_1\oplus S_2$, where $S_i\in\{{\mathfrak
sl}(2),\,W(1;\un{1}),\,H(2;\un{1})^{(2)}\}$ for $i=1,2$. Since
$L[\mu,\beta]\subset\Der \widetilde{S}$, the Lie algebra
$\pi(K(\mu,\beta))\supset\widetilde{S}$ is semisimple. Hence
$\ker\pi=\rad\,K(\mu,\beta)$. By Proposition~2.3(3), we have
$\overline{T}\subset\psi(H)$, which implies that no nonzero root
in $\Phi(L[\mu,\beta],T)$ vanishes on $\psi(H)$. It follows that
$K_{i\beta}\subset \ker\pi$ for all $i\in\F_p^*$. As
$\nu(K(\mu,\beta))$ is semisimple, we have
$\ker\pi=\rad\,K(\mu,\beta)\subset\ker\nu$. As
$TR(G[\widetilde{\mu}])=1$, either $(\nu\circ\pi^{-1})(S_1)=(0)$
or $(\nu\circ\pi^{-1})(S_2)=(0)$. No generality will be lost by
assuming that the latter case occurs. Then
$(\nu\circ\pi^{-1})(S_1)=\,U\cong S_1$.

Denote by ${\mathfrak m}$ the centralizer of $S_1$ in
$L[\mu,\beta]$. This is an ideal of $L[\mu,\beta]$ containing
$S_2$. Since ${\mathfrak m}/S_2$ is solvable, the preceding remark
implies that $\pi^{-1}({\mathfrak m})\subset\ker \nu$. When $S_2$
is Witt (resp., Hamiltonian), the subalgebra
$Q[\mu,\beta]+{\mathfrak m}$ has codimension $1$ (resp., $2$) in
$L[\mu,\beta]$. As in the previous case we now obtain that
$(\nu\circ\pi^{-1})\big(Q[\mu,\beta]\big)$ is a maximal
compositionally classical subalgebra of $G[\widetilde{\mu}]$. It
contains $\nu(K(\beta))$ because $\beta$ is solvable. Thus, all
$R$-roots of $G$ are proper. Now proceed as in the original proof
of \un{Lemma~5.2}.\qed

\smallskip

The rest of \un{Section~5} and most of \un{Section~6} are
essentially self-contained: they rely only on earlier results in
\cite{St93} and all arguments hold for $p>3$. However, in what
follows we will need a slightly different version of
\un{Theorem~6.7}.

According to \un{Theorem~3.5}, the maximal subalgebra $T+Q$ of
$T+L$ gives rise to a long standard filtration in $T+L$, and the
corresponding graded Lie algebra $G:={\rm gr}(T+L)$ has a unique
minimal ideal $A(L,T)$. Furthermore, the Lie algebra
$A(L,T)=\,\bigoplus_{i\in\Z}\,A_{[i]}$ is graded and there exist a
nonnegative integer $m$ and a simple graded Lie algebra
$S(L,T)=\,\bigoplus_{i\in\Z}\,S_{[i]}$ such that
$$A(L,T)=S(L,T)\ot\OO(m;\un{1}),\quad
A_{[i]}=S_{[i]}\ot\OO(m;\un{1})\quad\, (\forall\,i\in\Z)$$ as
graded Lie algebras. Our next result describes the graded
component $S_{[0]}$ of $S(L,T)$.
\begin{theo}
The Lie algebra $S_{[0]}$ is one of the following:
\begin{itemize}
\item[(a)] $1$-dimensional;

\smallskip

\item[(b)] classical simple;

\smallskip

\item[(c)] ${\mathfrak sl}(kp),\,$ ${\mathfrak gl}(kp)$ or
$\,{\mathfrak pgl}(kp)$ for some $k\ge 1$;

\smallskip

\item[(d)] $S'_{[0]}\oplus C$ where $C=C(S_{[0]})$ is
$1$-dimensional and $S'_{[0]}$ is either classical simple or
${\mathfrak pgl}(kp)$ for some $k\ge 1$.
\end{itemize}
\end{theo}
\pf (1) We know from \un{Proposition~3.10} and
\un{Proposition~6.5(4)} that $H_{[0]}={\mathfrak c}_{S_{[0]}}(T)$
is an abelian Cartan subalgebra of $S_{[0]}$. By
\un{Proposition~6.1(3)}, the semisimple quotients of the
$2$-sections of $S_{[0]}$ relative to $H_{[0]}$ are of types
$(0)$, ${\rm A}_1$, ${\rm A}_1\times {\rm A}_1$, ${\rm A}_2$,
${\rm C}_2$ or ${\rm G}_2$. In view of \un{Lemma~6.2(2)} and
\un{Proposition~6.5(1)} we have that
$\rad\,S_{[0]}(\bar{\alpha})\subset H_{[0]}$ for every root
$\bar{\alpha}\in\Phi(S_{[0]},H_{[0]})$. Consequently,
$S_{[0]}(\bar{\alpha})$ is nonsolvable if $\bar{\alpha}\ne 0$. By
\un{Proposition~3.10} and \un{Proposition~6.5(2)}, the torus $T$
acts on $S_{[-1]}$ and $H_{[0]}$ distinguishes the weight spaces
of $S_{[-1]}$ relative to $T$. Since those are $1$-dimensional, by
properties of $Q(L,T)$, we derive that every $S_{[0]}$-submodule
of $S_{[-1]}$ is $T$-stable. The construction of $S(L,T)$ now
yields that $S_{[-1]}$ is an irreducible $S_{[0]}$-module. But
then $\rad\,S_{[0]}\subset H_{[0]}$ acts on $S_{[-1]}$ by scalar
operators. As a result, $\rad\,S_{[0]}=\, C(S_{[0]}).$

Suppose $\bar{\alpha},\bar{\beta}\in\Phi(S_{[0]},H_{[0]})$ are
linearly independent. As $$\rad\, S_{[0]}(\bar{\gamma})= \{h\in
H_{[0]}\,|\,\,\bar{\gamma}(h)=0\}\qquad\,
\big(\forall\,\bar{\gamma}\in\Phi(S_{[0]},H_{[0]})\setminus\{0\}\big)$$
and $\ker\bar{\alpha}\cap\ker\bar{\beta}$ has codimension $2$ in
$H_{[0]}$, it must be that
$S_{[0]}[\bar{\alpha},\bar{\beta}]\not\cong{\mathfrak sl}(2)$. But
then $S_{[0]}[\bar{\alpha},\bar{\beta}]$ has type ${\rm A}_1\times
{\rm A}_1$, ${\rm C}_2$ or ${\rm G}_2$. This implies that
$\bar{\alpha}$ and $\bar{\beta}$ are linearly independent as
functions on $H_{[0]}\cap S_{[0]}^{(1)}$. As
$\dim\,S_{[0],\bar{\gamma}}=1$ for any nonzero
$\bar{\gamma}\in\Phi(S_{[0]},H_{[0]})$, it follows that every ideal
of $S_{[0]}^{(1)}$ is $H_{[0]}$-stable.

Since $S_{[0]}=S_{[0]}^{(1)}+H_{[0]}$ and $H_{[0]}$ is abelian,
the derived subalgebra $S_{[0]}^{(1)}$ is perfect. The preceding
remark implies that $\rad\,S_{[0]}^{(1)}=\,C(S_{[0]})\cap
S_{[0]}^{(1)}$. Put ${\mathfrak
g}:=S_{[0]}^{(1)}/\rad\,S_{[0]}^{(1)}$ and let ${\mathfrak h}$
denote the image of $H_{[0]}\cap S_{[0]}^{(1)}$ in $\mathfrak g$.
Obviously, $\mathfrak g$ is perfect and semisimple. The above
discussion shows that $\mathfrak h$ is an abelian Cartan
subalgebra of $\mathfrak g$ and the pair $({\mathfrak g},
{\mathfrak h})$ satisfies the Mills--Seligman axioms. Since $p>3$,
the main result of \cite{MS}  enables us to conclude that
$\mathfrak g$ is a direct sum of classical simple Lie algebras.

\smallskip

\noindent (2) \un{Theorem~3.3(4)} yields $Q_{(p-2)}\ne (0)$.
Furthermore, the proof of \un{Theorem~3.3} shows that
$Q_{(p-1)}\ne (0)$ provided that $L\ne Q_{(-1)}$. It also shows
that if $L=Q_{(-1)}$, then there exist root vectors $x\in
L\setminus Q$ and $u\in Q_{(1)}$ with
$\big[[x,y],Q_{(p-2)}\big]\ne (0)$. Since $G_{[-1]}\subset
A(L,T)$, it follows that $A_{[p-2]}\ne (0)$. As a consequence,
$S_{[3]}\ne (0)$. Now Lemmas~12.4.2--12.4.4 of \cite{BW} apply and
yield that $\mathfrak g$ is simple or zero. At this point we can
refer to \cite[Corollary~12.4.7]{BW} to complete the proof. (All
our references to \cite[Sect.~12]{BW} work for $p>3$; we are not
interested in the $p$-structure of $S_{[0]}$ which is also
discussed in \cite{BW}.) \qed

\smallskip

 Now we are ready to determine the Lie algebra $S(L,T)$.

\smallskip

\noindent \un{\bf Theorem~6.7} (new). {\it Let $L$ be a
finite-dimensional simple Lie algebra over $F$ and suppose that
all tori of maximal dimension in $L_p$ are standard. Let $T$ be an
optimal torus in $L_p$ and assume that $Q(L,T)\ne L$. Let $A(L,T)$
be the minimal ideal of the graded Lie algebra $G={\rm gr}(T+L)$,
and $S(L,T)$ the simple graded Lie algebra such that $A(L,T)\cong
S(L,T)\ot\OO(m;\un{1})$. Then $S(L,T)$ is a restricted simple Lie
algebra of Cartan type.}

\smallskip

\pf We will show that the conditions (a)-(d) of the Recognition
Theorem  apply to the graded Lie algebra $S(L,T)$; see
\cite[Theorem~0.1]{BGP} and \cite[Theorem~5.6.1]{St04}.

\smallskip

\noindent (a) Theorem~5.1 shows that the graded component
$S_{[0]}$ of $S(L,T)$ satisfies condition~(a) of the Recognition
Theorem.

\smallskip

\noindent (b) In part~(1) of the proof of Theorem~5.1 it was
explained that $S_{[-1]}$ is an irreducible $S_{[0]}$-module. This
means that condition~(b) of the Recognition Theorem holds for
$S(L,T)$.

\smallskip

\noindent (c) It follows from the definition of a standard
filtration that the Lie algebra $\bigoplus_{i<0}\,G_{[i]}$ is
generated by its subspace $G_{[-1]}$. Consequently, the Lie algebra
$\bigoplus_{i<0}\,S_{[i]}$ is generated by $S_{[-1]}$. If
$\big[x,S_{[-1]}\big]=(0)$ for some nonzero $x\in \bigoplus_{i\ge
0}\,S_{[i]}$, then the subspace $\bigoplus_{i\ge 0}\,S_{[i]}$
contains a nonzero ideal of $S(L,T)$. Since this contradicts the
simplicity of $S(L,T)$ we derive that condition~(c) of the
Recognition Theorem holds for $S(L,T)$.

\smallskip

\noindent (d) Now suppose that $\big[x,S_{[1]}\big]=(0)$ for some
nonzero $x\in S_{[-j]}$ with $j\ge 0$. Since the subspace
$$Y_{[-j]}:=\big\{y\in S_{[-j]}\,|\,\, \big[y,
S_{[1]}\big]=(0)\big\}$$ is $S_{[0]}$-stable, we may assume that
$x\in S_{{[-j]},\alpha}$ is a root vector for $T$.

\smallskip

\noindent Suppose $j=0$. Take any $h\in Y_{[0]}\cap H_{[0]}$ and
any $\gamma\in \Phi_{[-1]}(S,T)$. It follows from the definition
of $Q(L,T)$ and \un{Theorem~3.10} that there exist $u\in
S_{[-1],\gamma}$ and $v\in S_{[1],-\gamma}$ with $[u,v]\ne 0$. As
mentioned in the proof of Theorem~5.1, the abelian Lie algebra
$H_{[0]}$ acts on the  weight spaces of $S_{[-1]}$ relative to
$T$. Since all these weight spaces are $1$-dimensional and
$[h,v]=0$, we then have
$$\gamma(h)[u,v]=[[h,u],v]=-[[u,[h,v]]=0.$$
It follows that $h$ annihilates $S_{[-1]}$. But then $h=0$, forcing
$Y_{[0]}\cap H_{[0]}=(0)$. Combining this with \un{Proposition~6.1}
and the Engel--Jacobson theorem, we now deduce that the ideal
$Y_{[0]}$ of $S_{[0]}$ acts nilpotently on $S_{[-1]}$. The
irreducibility of $S_{[-1]}$ yields $Y_{[0]}=(0)$.

\smallskip

\noindent Suppose $j=1$. Since $Y_{[-1]}$ is an
$S_{[0]}$-submodule of $S_{[-1]}$ and
$\big[S_{[-1]},S_{[1]}\big]\ne (0)$ by part~(c) of this proof, the
irreducibility of $S_{[-1]}$ now yields $Y_{[-1]}=(0)$.

\smallskip

\noindent Suppose $j=2$. Then the $1$-section $L(\alpha)$ fits
into a $2$-section $L(\alpha,\beta)$ whose semisimple quotient is
isomorphic to $K(3;\un{1})$; see \un{Proposition~3.4}. Since
$S_{[-2],\alpha}=Fx$, it follows from \un{Lemma~3.1(1)} that that
$\big[x, G_{[1]}\big]\ne (0)$. Because
$G_{[1]}=S_{[1]}\ot\OO(m;\un{1})$ and $x$ identifies with an
element in $S_{[-2]}\ot F\subset S_{[-2]}\ot\OO(m;\un{1})$, we now
obtain that $\big[x, S_{[1]}\big]\ne (0)$, a contradiction. Hence
$Y_{[-2]}=(0)$. As $S_{[-j]}=(0)$ for $j>2$, by
\un{Proposition~3.4}, we have proved that all conditions of the
Recognition Theorem are satisfied for $S(L,T)$.

\smallskip

\noindent (e) Applying the Recognition Theorem we obtain that
$S(L,T)$ is either classical or of Cartan type or a Melikian Lie
algebra. As explained in part~(2) of the proof of Theorem~5.1, we
have that $S_{[3]}\ne (0)$. So $S(L,T)$ has an unbalanced grading,
hence cannot be classical. The natural grading of any Melikian
algebra has depth $3$ and height $>3$. As $S_{[-3]}=(0)$, it
follows that $S(L,T)$ is not of Melikian type. We conclude that
the graded Lie algebra $S(L,T)$ is isomorphic to a Cartan type Lie
algebra $X(r;\un{s})^{(2)}$ regarded with its natural grading
(here $X\in\{W,S,H,K\}$ and $\un{s}=(s_1,\ldots,s_r)\in{\mathbb
N}^r$).

To show that $S(L,T)$ is restricted we  take any root vector $x\in
S_{[i],\alpha}$ with $i\in\{-1,-2\}$ and let $\beta$ be any $T$-root
of $S=S(L,T)$. The semisimple quotients of the $2$-sections of $S$
relative to $T$ are described in \un{Theorem~3.11}. Furthermore, the
proof of \un{Theorem~3.11} in conjunction with Theorem~5.1 shows
that $2$-sections of type~(7) do not occur. Since $\alpha\in\Phi_-$,
it also follows from the proof of \un{Theorem~3.11} that
$$(\ad x)^p(S(\alpha,\beta))\subset \rad_T\, S(\alpha,\beta)\subset
\textstyle{\bigoplus}_{i\ge 0}\,S_{[i]}(\alpha,\beta).$$ But then
$(\ad x)^p$ maps $S(L,T)$ into $\textstyle{\bigoplus}_{i\ge
0}\,S_{[i]}.$ Since $S(L,T)\cong X(r;\un{s})^{(2)}$ as graded Lie
algebras, this forces $\un{s}=\un{1}$. As a consequence, $S(L,T)$ is
restricted; see \cite[Corollary~7.2.3]{St04} for example.\qed

\section{Classification results}

Similar to \cite[Sect.~7]{St93} the determination of the Lie
algebra $S(L,T)$ allows one to classify a large family of
finite-dimensional simple Lie algebras. Note that all results and
arguments used in \cite[Sect.~7]{St93} are valid for $p>3$.
\begin{theo}{\rm (cf.~\cite[Theorem~7.3]{St93}).} Let $L$ be a
finite-dimensional simple Lie algebra over $F$ such that all tori
of maximal dimension in $L_p$ are standard. Let $T$ be an optimal
torus in $L_p$ and suppose that $Q(L,T)\ne L$ and
$L[\alpha,\beta]^{(1)}\not\cong H(2;\un{1};\un\Phi(\tau))^{(1)}$
for any two roots $\alpha,\beta\in\Gamma(L,T)$. Then $L$ is
isomorphic to a Cartan type Lie algebra and $Q(L,T)$ is contained
in the standard maximal subalgebra of $L$.
\end{theo}
\pf One argues as in the original proof of Theorem~7.3 in
\cite{St93} to construct a maximal subalgebra $L_{(0)}$ containing
$Q(L,T)$ and to show that the pair $(L,L_{(0)})$ satisfies all
conditions of the Recognition Theorem for filtered Lie algebras;
see \cite[Theorem~5.6.2]{St04}. The argument in \cite[pp.~57,
58]{St93} shows that $L=L_{(-2)}$. This implies that $L$ is not
isomorphic to a Melikian algebra. Since $S(L,T)$ is of Cartan
type, it follows from the construction of $L_{(0)}$ that
$L_{(p-2)}\ne (0)$. But then $L$ cannot be classical. By the
Recognition Theorem, $L$ must be isomorphic as a filtered Lie
algebra to a Cartan type Lie algebra regarded with its standard
filtration. This completes the proof. \qed

\smallskip

We continue assuming that all tori of maximal dimension in $L$ are
standard and $Q(L,T)\ne L$. In view of Theorem~6.1 we can also
assume now that for any optimal torus $T$ in $L_p$ there are
$\alpha,\beta\in\Gamma(L,T)$ such that $L[\alpha,\beta]^{(1)}\cong
H(2;\un{1};\Phi(\tau))^{(1)}$. For $p>7$, the simple Lie algebras
with these properties are classified in \cite{St94}. We will go
through the arguments in \cite{St94} to verify whether they are
still valid for $p=5,7$. All our references to \cite{St94} will be
underlined.

We have already shown in \cite{PS4} that the results of
\cite{St91a} hold for $p=5,7$. Note that \cite{St91a} is the main
prerequisite to \cite{St94}. Inspection shows that all results and
arguments used in \cite{St94} are valid for $p>5$. In fact, only
one minor issue in \cite[Sect.~2]{St94} requires our attention; it
arises when $p=5$.

\smallskip
\noindent \un{\bf Proposition~2.4} (new parts).

\begin{itemize}
\item[(3a)]\, {\it If $\alpha\in\Phi_{[-1]},\,$
$\beta\in\Phi_{[0]}$ and $\alpha+\beta\in\Phi_{[0]}$, then
$\alpha$ is Hamiltonian, $p=5$, and $\,\beta=2\alpha$.}

\smallskip

\item[(3d)]\, $[Q_\mu,L_\lambda]\cap Q\subset Q_{(1)}$ {\it for
all } $\lambda\in\Phi_{[-1]}$ and $\mu\in\Delta$.
\end{itemize}

\pf (3a) The assertion follows from \cite[Lemma~5.5.1]{St04}.

\smallskip

\noindent (3d) Recall that $\Phi_{[-1]}({\mathrm
gr}(T+L),T)=\Phi_{[-1]}(S,T)+\Delta$; see \un{Proposition~2.4(2)}.
Since $\mu+\lambda\in\Phi_{[-1]}+\Delta=\Phi_{[-1]}$ and
$\Phi_{[-1]}\cap \Phi_{[0]}=\emptyset$ by
\un{Proposition~2.4(3b)}, the statement follows.\qed

\medskip

As a result of the above changes we have to modify slightly the
statement and the original proof of \un{Proposition~2.5}: the
subspace $V_0$ from \un{Proposition~2.5} has to be selected in a
more sophisticated fashion.

\smallskip

\noindent \un{\bf Proposition~2.5} (new proof).
\begin{itemize}
\item[1)] {\it There exists a $T$-invariant subspace
$V_{-1}\subset L$ such that
$$L=V_{-1}+Q,\quad V_{-1}\cap Q=(0).$$} \item[2)] {\it There
exists a $T$-invariant subspace $V_0\subset Q$ such that}
$$\quad Q_{(1)}\cap V_0=(0)\ \mbox{ {\it and} }\ \,
(V_0+Q_{(1)})/Q_{(1)}=\,A(L,T)_{[0]}.$$

\item[3)] {\it For $i=-1,0$, let $R_i\subset V_i$ denote the
preimage of $S_{[i]}\ot\OO(m;\un{1})_{(1)}$ under the linear
isomorphism $V_i\,\stackrel{\sim}{\rightarrow}\,
(V_i+Q_{(i+1)})/Q_{(i+1)}= A(L,T)_{[i]}$. Then the following
statements hold:}

\medskip

\begin{itemize}
\item[a)] $V_{-1}\subset \sum_{\mu\in\Phi_{[-1]}}L_\mu;$

\smallskip

\item[b)] $\sum_{\mu\not\in\Delta}Q_\mu\subset V_0+Q_{(1)}$;

\smallskip

\item[c)] $Q=V_0+Q(\Delta)+Q_{(1)}$ {\it and }  $V_0\subset
\sum_{\mu\in\Phi_{[0]}}Q_\mu+Q_{(1)}$;

\smallskip

\item[d)] $[V_{-1},V_{-1}]\subset Q$;

\smallskip

\item[e)] $[T+V_0+Q(\Delta), V_{-1}]\subset V_{-1}+Q_{(1)}$;

\smallskip

\item[f)] $[Q_{(1)}, V_{-1}]\subset V_0+Q_{(1)}$;

\smallskip

\item[g)] $[T+Q,V_0]\subset V_0+Q_{(1)}$;

\smallskip

\item[h)] $[V_0,R_{-1}]\subset R_{-1}+Q_{(1)}\subset
[V_0,R_{-1}]+Q_{(1)}$;

\smallskip

\item[i)] $[V_0,R_0]\subset R_0+Q_{(1)}\subset [V_0,R_0]+Q_{(1)}$;

\smallskip

\item[j)] $V_0+Q_{(1)}$ {\it is an ideal of} $\,Q$;

\smallskip

\item[k)] $[R_0,V_{-1}]\subset R_{-1}+Q_{(1)}\subset
[R_0,V_{-1}]+Q_{(1)}.$
\end{itemize}
\end{itemize}

\pf The original proof goes through for $p>5$.  So we assume from
now that $p=5$. We choose $V_{-1}$ as in the original proof. Then
assertions~1), 3a) and 3d) hold.  Let $\Phi'_{[0]}$ denote the set
of all nonzero $T$-roots of $A(L,T)_{[0]}$. For every
$\mu\in\Phi'_{[0]}$ choose a nonzero $u_\mu\in Q_\mu$ such that
$A(L,T)_{[0],\mu}=\,F\bar{u}_\mu$, where $\bar{u}_\mu$ stands for
the coset of $u_\mu$. If $\mu$ is not Hamiltonian, set
$v_\mu:=u_\mu$. If $\mu$ is Hamiltonian, then $\pm 3\mu$ are
$T$-weights of $L/Q$, so that $V_{-1,\pm 3\mu}\ne (0)$. Pick $u_{\pm
3\mu}\in V_{-1,\pm 3\mu}\setminus\{0\}$. As $p=5$, we have
$[u_{3\mu},u_\mu]\in L_{-\mu}$. Then
$[u_{3\mu},u_\mu]=ru_{-\mu}+q_{-\mu}$ for some $r\in F$ and
$q_{-\mu}\in Q_{(1),-\mu}$. Since $\rad\,L(\mu)\subset Q_{(1)}$ and
the image of $L(\mu)\cap Q_{(1)}$ in $L[\mu]\subset
H(2;\un{1})^{(1)}$ contains ${H(2;\un{1})^{(2)}}_{(1)},$ it is easy
to see that there is $w_{\mu}\in Q_{(1),\mu}$ such that
$$[u_{3\mu},w_{\mu}]\equiv\, ru_{-\mu}\ \,({\mathrm mod}\,\,
Q_{(1)}).$$ Put $v_\mu:=u_\mu-w_{\mu}$. Then $[u_{3\mu}, v_\mu]\in
Q_{(1)}.$ Now set
$$V_0:=\,(H_{[0]}\ot F)\oplus\,
\textstyle{\bigoplus}_{\mu\in\Phi'_{[0]}}\, Fv_\mu.$$ By
construction, we have $Q_{(1)}\cap V_0=(0)$ and
$(V_0+Q_{(1)})/Q_{(1)}=\,A(L,T)_{[0]}$. As $Q_{(1)}\cap V_0=(0)$,
this yields assertion~2). In view of the new
\un{Proposition~2.4(3a)} our choice of $V_0$ ensures that
$$[V_0,V_{-1}]\subset V_{-1}+Q_{(1)}.$$ To prove assertions~3b), 3c),
3f), 3g) and 3j) one can argue as in the original proof.
Assertion~3e) follows from the new \un{Propostion~2.4(3d)} and the
displayed inclusion. Assertions 3h), 3i) and 3k) follows from the
displayed inclusion and 3e) by the same argument as in the
original proof. \qed

\medskip

\noindent \un{\bf Lemma~2.6} (new). {\it Let $u\in L_{\alpha}$,
$f\in (V_0+Q(\Delta))_{\beta}$ and $v\in L_{\alpha+\beta}$ be such
that $u\not\in Q$ and $[f,u]-v\in Q$. Then $[f,u]-v\in Q_{(1)}$.}
\pf Write $u=u_\alpha+u'$ with $u_\alpha\in V_{-1}$ and $u'\in
Q_\alpha$. Then $[f,u_\alpha]\in V_{-1}+Q_{(1)}$ thanks to
\un{Proposition~2.5(3e)}, while \un{Proposition~2.4(3b)} yields
$u'\in Q_{(1)}$. As $Q_{(1)}$ is an ideal of $Q$ and
$V_0+Q(\Delta)\subset Q$, we have $[f,u']\in Q_{(1)}$. So it remains
to show that $[f,u_\alpha]-v\in Q_{(1)}$.

If $v\not\in Q$, then the coset of $v$ spans $
L_{\alpha+\beta}/L_{\alpha+\beta}\cap Q_{(1)}$, again by
\un{Proposition~2.4(3b)}. Hence the assertion holds in this case.
Now suppose $v\in Q$. If $v\in Q_{(1)}$, we are done; so suppose for
a contradiction that $v\not\in Q_{(1)}$. Then the new
\un{Proposition~2.4(3d)} yields that $\alpha$ is Hamiltonian, $p=5$,
and $\beta=2\alpha$. On the other hand, it is immediate from
\un{Corollary~2.3(3)} that $2\Phi_{[-1]}\cap \Delta=\emptyset$. This
forces $f\in V_0$. But then
$$[f,u_\alpha]\in [V_{-1},V_0]\cap Q\subset Q_{(1)},$$
by our choice of $V_0$. This completes the proof. \qed

\medskip

With these substitutions, all arguments in \cite{St94} work. As a
result, we obtain

\smallskip
\noindent \un{\bf Theorem~6.4} (new). {\it Let $L$ be a
finite-dimensional simple Lie algebra over $F$ such that all tori
of maximal dimension in $L_p$ are standard. Let $T$ be an optimal
torus in $L_p$ and suppose that $Q(L,T)\ne L$ and
$L[\alpha,\beta]^{(1)}\cong H(2;\un{1};\Phi(\tau))^{(1)}$ for some
$\alpha,\beta\in\Gamma(L,T)$. Then $L$ is isomorphic to a Lie
algebra of Cartan type.}

\medskip

We summarize the results of this paper and of \cite{PS4} as
follows:
\begin{theo}
Let $L$ be a finite-dimensional simple Lie algebra over an
algebraically closed field of characteristic $p>3$. If all tori of
maximal dimension in the the semisimple $p$-envelope $L_p$ of $L$
are standard, then $L$ is isomorphic to a Lie algebra of classical
 or Cartan type.
\end{theo}
\pf If $Q(L,T)=L$, then all roots in $\Gamma(L,T)$ are either
solvable or classical. So the assertion follows from
\cite[Theorems~C and D]{PS4} in this case. If $Q(L,T)\ne L$, the
assertion follows from Theorem~6.1 and the new
\un{Theorem~6.4}.\qed

\begin{rem}
{\rm The argument in the last two paragraphs of \cite[p.~792]{PS4}
can be streamlined as follows: When $\alpha$ is solvable with
$\alpha(H)\ne 0$, Proposition~3.8 of \cite{PS4} yields that
$[L_{\alpha},L_{-\alpha}]$ consists of $p$-nilpotent elements of
$L_p$. But then $[L_{\alpha},L_{-\alpha}]\subset {\rm
nil}\,H\subset{\rm nil}\,\widetilde{H}$, contrary to our choice of
$\alpha$. Therefore, $[L_{\gamma},L_{-\gamma}]\subset{\rm
nil}\,\widetilde{H}$ whenever $\gamma(\widetilde{H}^{(1)})\ne 0$.}
\end{rem}

\end{document}